\documentclass[12pt,reqno]{amsart}

\usepackage[T1]{fontenc}
\usepackage{lmodern}
\usepackage{microtype}
\usepackage{mathrsfs}
\usepackage{amsmath,amssymb,amsfonts,mathtools}
\usepackage{booktabs,array,enumitem}
\usepackage{tikz-cd}
\usepackage[colorlinks=true,linkcolor=blue!55!black,citecolor=blue!55!black,urlcolor=blue!55!black]{hyperref}
\usepackage{aliascnt}
\usepackage[nameinlink,noabbrev]{cleveref}
\usepackage[a4paper,margin=30mm]{geometry}

\setlist[itemize]{leftmargin=2.0em,itemsep=2pt,topsep=3pt}
\setlist[enumerate]{leftmargin=2.2em,itemsep=2pt,topsep=3pt}
\numberwithin{equation}{section}
\allowdisplaybreaks
\newtheorem{theorem}{Theorem}[section]
\newaliascnt{proposition}{theorem}
\newtheorem{proposition}[proposition]{Proposition}
\aliascntresetthe{proposition}
\newaliascnt{lemma}{theorem}
\newtheorem{lemma}[lemma]{Lemma}
\aliascntresetthe{lemma}
\newaliascnt{corollary}{theorem}
\newtheorem{corollary}[corollary]{Corollary}
\aliascntresetthe{corollary}
\newaliascnt{claim}{theorem}

\aliascntresetthe{claim}
\theoremstyle{definition}
\newaliascnt{definition}{theorem}

\aliascntresetthe{definition}
\newaliascnt{construction}{theorem}

\aliascntresetthe{construction}
\theoremstyle{remark}
\newaliascnt{remark}{theorem}
\newtheorem{remark}[remark]{Remark}
\aliascntresetthe{remark}

\newcommand{\R}{\mathbb R}
\newcommand{\Z}{\mathbb Z}

\newcommand{\Hh}{\mathbb H}
\newcommand{\cA}{\mathcal A}
\newcommand{\cD}{\mathcal D}
\newcommand{\cE}{\mathcal E}
\newcommand{\cG}{\mathcal G}
\newcommand{\cT}{\mathcal T}
\newcommand{\cH}{\mathcal H}
\newcommand{\cK}{\mathcal K}
\newcommand{\cQ}{\mathcal Q}
\newcommand{\cN}{\mathcal N}
\newcommand{\mathsfH}{\mathsf H}
\newcommand{\Diff}{\operatorname{Diff}}
\newcommand{\Map}{\operatorname{Map}}

\newcommand{\Ad}{\operatorname{Ad}}
\newcommand{\im}{\operatorname{Im}}
\newcommand{\Coker}{\operatorname{Coker}}
\newcommand{\Hess}{\operatorname{Hess}}
\newcommand{\Scal}{\operatorname{Scal}}

\newcommand{\Sec}{\operatorname{sec}}
\newcommand{\op}{\mathrm{op}}
\newcommand{\id}{\mathrm{id}}
\newcommand{\supp}{\operatorname{supp}}
\newcommand{\Int}{\operatorname{Int}}
\newcommand{\pr}{\operatorname{pr}}
\newcommand{\Span}{\operatorname{span}}
\newcommand{\dd}{\,\mathrm d}
\newcommand{\eps}{\varepsilon}
\newcommand{\norm}[1]{\left\lVert #1\right\rVert}
\newcommand{\abs}[1]{\left\lvert #1\right\rvert}
\newcommand{\ip}[2]{\left\langle #1,#2\right\rangle}

\title[Hitchin's invariant and homotopy 9-spheres]{Hitchin's $\alpha$-invariant and positive sectional curvature on homotopy $9$-spheres}

\author{Yang-Hui He}
\address{London Institute for Mathematical Sciences, Royal Institution, London W1S 4BS, United Kingdom}
\address{Merton College, University of Oxford, Oxford OX1 4JD, United Kingdom}
\email{yh@lims.ac.uk}

\author{Ziran Liu}
\address{Shanghai Institute for Mathematics and Interdisciplinary Sciences (SIMIS), Shanghai 200433, China}
\address{Research Institute of Intelligent Complex Systems, Fudan University, Shanghai 200433, China}
\email{zliu@simis.cn}

\author{Shing-Tung Yau}
\address{Yau Mathematical Sciences Center, Tsinghua University, Beijing 100084, China}
\email{styau@mail.tsinghua.edu.cn}


\hypersetup{
  pdftitle={Hitchin's alpha-invariant and positive sectional curvature on homotopy 9-spheres},
  pdfauthor={Yang-Hui He, Ziran Liu, and Shing-Tung Yau},
  pdfsubject={Positive sectional curvature and the classification of smooth homotopy 9-spheres},
  pdfkeywords={positive sectional curvature, exotic sphere, homotopy sphere, Kervaire sphere, Gromoll filtration, mapping class group}
}

\begin{document}

\begin{abstract}
We prove that a smooth homotopy $9$-sphere admits a Riemannian metric of
positive sectional curvature if and only if its Hitchin $\alpha$-invariant
vanishes, or equivalently if and only if it admits a Riemannian metric of
positive scalar curvature.
\end{abstract}

\maketitle
\enlargethispage{2pt}

\section{Introduction}

For a closed spin $n$-manifold $M$, Hitchin's $\alpha$-invariant is the
real index of the spin Dirac operator, taking values in the real
$K$-homology group of a point:
\[
 \alpha(M)\in KO_n(\mathrm{pt}).
\]
Hitchin proved that this real Dirac index vanishes whenever $M$ admits
a metric of positive scalar curvature \cite{Hitchin}.  Conversely,
Stolz proved that its vanishing is sufficient when $M$ is also simply
connected and has dimension at least five \cite{Stolz}.

Positive sectional curvature implies positive scalar curvature.  For
smooth homotopy $9$-spheres, we prove that the converse holds at the
level of existence of metrics.  A smooth homotopy $9$-sphere is a closed
smooth manifold homotopy equivalent to $S^9$; it is \emph{exotic} if it
is not diffeomorphic to the standard sphere.
Every homotopy $9$-sphere is orientable; after an orientation is fixed,
its spin structure is unique.
Dimension nine is the first dimension in which a simply connected
homotopy sphere can have nonzero $\alpha$-invariant
\cite[Section 2.3]{CrowleySchick}.

For a Riemannian metric $g$, write $\langle\cdot,\cdot\rangle_g$ for
its tangent-space inner product and let $\nabla$ be its Levi--Civita
connection.  We use the curvature convention
\[
 R^g(X,Y)Z=\nabla_X\nabla_YZ-\nabla_Y\nabla_XZ-\nabla_{[X,Y]}Z.
\]
For linearly independent $X,Y\in T_pM$, the sectional curvature of
$\Pi=\operatorname{span}\{X,Y\}$ is
\begin{equation}\label{eq:intro-sectional}
 \Sec_g(\Pi)=
 \frac{\langle R^g(X,Y)Y,X\rangle_g}
  {\langle X,X\rangle_g\langle Y,Y\rangle_g-\langle X,Y\rangle_g^2}.
\end{equation}
The notation $\Sec_g>0$ means that this quantity is positive at every
point and on every two-plane.  We write $\Scal_g$ for scalar curvature;
for a $g$-orthonormal frame $(e_1,\ldots,e_n)$ it satisfies
$\Scal_g=2\sum_{i<j}\Sec_g(\operatorname{span}\{e_i,e_j\})$.
All metrics considered below are smooth.

\subsection{Main result}

\begin{theorem}\label{thm:main-intro}
Let $\Sigma$ be a smooth homotopy $9$-sphere.  The following conditions
are equivalent:
\begin{enumerate}[label=\textup{(\roman*)}]
 \item $\Sigma$ admits a Riemannian metric of positive sectional curvature;
 \item $\Sigma$ admits a Riemannian metric of nonnegative sectional curvature;
 \item $\Sigma$ admits a Riemannian metric of positive scalar curvature;
 \item $\alpha(\Sigma)=0\in KO_9(\mathrm{pt})\cong\Z/2$.
\end{enumerate}
\end{theorem}

To prove \textup{(iv)}$\Rightarrow$\textup{(i)}, we construct a positively
curved metric for each nonstandard smooth type with vanishing index.  The standard sphere has its round metric.
The implication \textup{(iii)}$\Rightarrow$\textup{(iv)} is Hitchin's
obstruction.  Its extension to nonnegative sectional curvature is
proved in \cref{prop:no-nonnegative}.

To list the smooth types, let $\Theta_9$ be the group of oriented
homotopy $9$-spheres under connected sum.  The classical classification
gives
\begin{equation}\label{eq:intro-theta}
 \Theta_9\cong(\Z/2)^3.
\end{equation}
Write $\pi_j^s=\varinjlim_k\pi_{j+k}(S^k)$ for the $j$th stable
homotopy group of spheres, with transition maps given by suspension.
Let $SO=\varinjlim_k SO(k)$, where the inclusions add an identity
coordinate, and let $J_j:\pi_j(SO)\to\pi_j^s$ be the stable
$J$-homomorphism.  We put $\Coker J_j=\pi_j^s/\im J_j$ and use
$[u]$ for the image of $u\in\pi_j^s$ in this quotient.  The subgroup $bP_{10}$ consists of
homotopy $9$-spheres that bound parallelizable $10$-manifolds.  There
is a split exact sequence
\begin{equation}\label{eq:intro-split}
 0\longrightarrow bP_{10}\longrightarrow\Theta_9
 \xrightarrow{\,p_9\,}\Coker J_9\longrightarrow0,
\end{equation}
with
\begin{equation}\label{eq:intro-stems}
 bP_{10}=\Z/2\{\kappa\},\qquad
 \Coker J_9=\Z/2\{[\eta\epsilon]\}\oplus\Z/2\{[\mu]\}.
\end{equation}
Here $\kappa$ denotes the oriented class of the Kervaire sphere,
the boundary of the plumbing of two tangent disc bundles of $S^5$.  We use Toda's notation
$\eta\in\pi_1^s$, $\nu\in\pi_3^s$, $\epsilon\in\pi_8^s$, and
$\mu\in\pi_9^s$.  The relation $[\eta\epsilon]=[\nu^3]$ holds in
$\Coker J_9$; the relevant stable homotopy calculations are recalled in
Section~\ref{sec:classification} \cite{KM,Toda,Wang}.

Our construction specifies a lift $\lambda\in\Theta_9$ of
$[\eta\epsilon]$.  Choose $\widehat\mu\in\Theta_9$ with
$p_9(\widehat\mu)=[\mu]$.  Then
\begin{equation}\label{eq:intro-decomp}
 \Theta_9=\langle\kappa\rangle\oplus\langle\lambda\rangle
              \oplus\langle\widehat\mu\rangle,\qquad
 \ker\alpha=\{0,\kappa,\lambda,\kappa+\lambda\}.
\end{equation}
The four classes in $\widehat\mu+\ker\alpha$ have nonzero index
\cite{Hitchin,Wang,Sarkar}.  Table~\ref{tab:classification} records
\cref{thm:main-intro} for all eight classes.  Orientation reversal does
not change these entries, since every element of $\Theta_9$ has order
at most two.

\begin{table}[htbp]
\centering\small
\renewcommand{\arraystretch}{1.16}
\begin{tabular}{@{}c c c c@{}}
\toprule
Class in $\Theta_9$ & Image under $p_9$ & $\alpha$ &
 \begin{tabular}[c]{@{}c@{}}Admits positive\\sectional curvature?\end{tabular}\\
\midrule
$0$ & $0$ & $0$ & Yes\\
$\kappa$ & $0$ & $0$ & Yes\\
$\lambda$ & $[\eta\epsilon]$ & $0$ & Yes\\
$\kappa+\lambda$ & $[\eta\epsilon]$ & $0$ & Yes\\
\midrule
$\widehat\mu$ & $[\mu]$ & $1$ & No\\
$\widehat\mu+\kappa$ & $[\mu]$ & $1$ & No\\
$\widehat\mu+\lambda$ & $[\mu]+[\eta\epsilon]$ & $1$ & No\\
$\widehat\mu+\kappa+\lambda$ & $[\mu]+[\eta\epsilon]$ & $1$ & No\\
\bottomrule
\end{tabular}
\caption{The classification of positive sectional curvature on the eight
oriented smooth homotopy $9$-spheres.}
\label{tab:classification}
\end{table}

\subsection{Related work}
Chen, Ni, and Zhang \cite{CNZ} construct a positively curved metric on the Kervaire $9$-sphere.  Subsequently, Chen, Deng, Ni, and Zhang \cite[Theorem 1.2]{CDNZ} extended their preprint to obtain the same classification as \cref{thm:main-intro}.  

A complete draft of the present work had already been prepared before that revision appeared. The two approaches share certain topological ingredients, including
a four-periodic normal-framing calculation, but differ in their geometric constructions.  In particular, the treatment of the
remaining nonstandard class with vanishing $\alpha$-invariant in \cite[Sections 10--12]{CDNZ} requires a positively curved cylinder,
a four-form modification, and two disc fillings.  In contrast,
our invariant-transition theorem \cref{thm:invariant-transition}
provides a direct two-disc construction, yielding a unified geometric
proof for all three nonstandard classes in $\ker\alpha$.

The argument here realizes all three nonstandard classes by gluing two
quotient discs directly.  We identify the Kervaire class through
Speran\c ca's pullback model \cite{Speranca}, and identify $\lambda$
using a family defining his exotic $8$-sphere.  The third attaching map
comes from left quaternionic multiplication and an invariant transition
function.  It is treated by \cref{thm:invariant-transition}.

The metric construction uses the two-disc method of He--Liu--Yau
\cite{HLY}.  Galaz-Garc\'ia \cite{GG} and Deng--Hu--Zhang
\cite{DHZ} give versions for conjugation-equivariant transition
functions; the latter allows arbitrary orthogonal $S^3$-representations.  These supply the geometric framework for the
classes $\kappa$ and $\lambda$.  We give the curvature estimates and boundary comparisons for both types of transition in Section~\ref{sec:geometry}.

\subsection{Outline of the proof}
We first construct three attaching diffeomorphisms of $S^8$ and identify
the smooth homotopy spheres obtained by gluing two $9$-discs along them.
The maps are defined by left and right quaternionic multiplication in a
tubular neighborhood of a standard $S^3\subset S^8$.  Their formulas,
support conditions and normal coordinates are fixed in
Section~\ref{sec:tubular-neighborhood}.

The first sphere is identified with Speran\c ca's Kervaire model.  The
second is constructed from a family of disc diffeomorphisms defining
his exotic $8$-sphere.  Precomposition with the Hopf map on the
classifying space of the diffeomorphism group gives a $9$-sphere with
image $[\eta\epsilon]$ under $p_9$.  We denote this class by $\lambda$.
A relative mapping-class calculation and quaternionic inversion
identify the third sphere as $\kappa+\lambda$.

For each attaching map we construct positively curved metrics on the
two discs.  The southern metric is a quotient of a connection metric,
and the northern metric is a quotient of a doubly warped product.
For the southern connection metric over $(B,g_B)$, the mixed curvature
contains $-r^{-1}\Hess_{g_B}r$, where $r$ is the fiber radius.
A strongly concave function $\varphi:B\to\R$ supplies the positive
term when $r^2=\eps e^{\eps\varphi}$; the scale $\eps>0$ is then
chosen small relative to the connection curvature and its covariant
derivative.
Both quotient actions are free.  The estimates apply to every
horizontal two-plane, including those over points where the base action
has a nontrivial stabilizer.

Let $\phi:\partial D_N\to\partial D_S$ be the attaching map already
used to identify the smooth sphere.  If $h_N,h_S$ are the induced
boundary metrics and $B_N,B_S$ their outward second fundamental forms,
the constructed metrics satisfy
\[
 h_N=\phi^*h_S,\qquad B_N+\phi^*B_S>0.
\]
Here $B(Y,Z)=\langle\nabla_Yn,Z\rangle_g$ for the outward unit normal $n$, and
$>0$ means positive definiteness.  The Reiser--Wraith theorem
\cite[Theorem A(i), $k=1$]{RW} then gives a smooth positively curved
metric on the specified sphere.  Section~\ref{sec:classification} fixes
the disc parametrizations used in these boundary comparisons.

Sections~\ref{sec:classes-kappa-lambda} and \ref{sec:last-topology} identify the three
smooth types; Section~\ref{sec:geometry} constructs the disc metrics, and
Section~\ref{sec:applications} completes the proof.  The appendices
contain the connection-metric estimates, warping functions, and
parameterized smoothing comparison.

\subsection*{Acknowledgments}
The authors used ChatGPT lightly to assist in exploration and reviewing mathematical arguments and in preparing the manuscript. The authors take full responsibility for the mathematical content and the final text.

\section[Homotopy 9-spheres and clutching conventions]{\texorpdfstring{Homotopy $9$-spheres and clutching conventions}{Homotopy 9-spheres and clutching conventions}}\label{sec:classification}

\subsection[The group of homotopy 9-spheres]{\texorpdfstring{The group $\Theta_9$}{The group Theta 9}}
For the dimensions $n=8,9$ used below, let $\Theta_n$ denote the group of
oriented diffeomorphism classes of smooth homotopy $n$-spheres, with
connected sum as group operation.  The clutching construction identifies
this group with orientation-preserving diffeomorphisms of $S^{n-1}$ modulo
those extending over $D^n$ \cite{KM,Cerf}.  We use the convention
\[
   \Sigma_\varphi=D^n_N\cup_\varphi(-D^n_S):
\]
$\varphi$ is orientation preserving in the two standard boundary
parameters, while the southern disc is given the opposite orientation.
With this convention composition represents connected sum.  We suppress
the orientation sign in subsequent gluing expressions.  Reversing the
convention negates the sphere class, which leaves it unchanged in
$\Theta_8$ and $\Theta_9$.

The Kervaire--Milnor map $p_n$ records a stable normal framing modulo
the stable $J$-image.  More precisely, a homotopy $n$-sphere admits a
stable normal framing.  The Pontryagin--Thom collapse associated with
such a framing represents an element of $\pi_n^s$; changing the framing
changes this element by an element of $\im J_n$.  Its class in
$\Coker J_n$ is $p_n([\Sigma])$ \cite{KM}.  We use this convention for
$n=8,9$ and write $p_n(\Sigma)$ for $p_n([\Sigma])$.

The group structure and exact sequence are given in
\eqref{eq:intro-theta}--\eqref{eq:intro-stems}.
We use the stable-stem notation of Toda and Wang; in particular,
$\bar\nu\in\pi_8^s$.  In degree nine,
\begin{equation}\label{eq:pi9stable}
 \pi_9^s=\Z/2\{\eta\bar\nu\}\oplus
         \Z/2\{\eta\epsilon\}\oplus
         \Z/2\{\mu\},
 \qquad
 \im J_9=\Z/2\{\eta\bar\nu+\eta\epsilon\}.
\end{equation}
The degree-eight and degree-nine entries of
\cite[Appendix C, Table 4]{Wang} give the displayed groups and
$\im J_9$.  Relation C.2(1) there is $\eta\bar\nu=\nu^3$;
hence $[\eta\epsilon]=[\nu^3]$ in $\Coker J_9$.

In degree eight, $\Theta_8\cong\Z/2$, and the Kervaire--Milnor map
$p_8:\Theta_8\to\Coker J_8\cong\Z/2\{[\epsilon]\}$ is an
isomorphism \cite{KM,Wang}.  We write $\Sigma^8_{[\epsilon]}$ for its
nonzero class.
The two lifts of $[\eta\epsilon]$ differ by $\kappa$;
Section~\ref{sec:classes-kappa-lambda} specifies $\lambda$ by a geometric construction.

\subsection[The alpha-invariant]{\texorpdfstring{The $\alpha$-invariant}{The alpha-invariant}}
Recall that the Hitchin invariant $\alpha(M)\in KO_j(\mathrm{pt})$
of a closed spin $j$-manifold is the real Dirac index.  It depends only
on the spin bordism class.  We abbreviate $KO_j(\mathrm{pt})$ to $KO_j$.
A homotopy $9$-sphere is orientable and has vanishing second
Stiefel--Whitney class because $H^2(\Sigma;\Z/2)=0$.  It is therefore
spin after an orientation has been fixed, and the corresponding spin
structure is unique because $H^1(\Sigma;\Z/2)=0$.  The real Dirac index defines a homomorphism
\[
  \alpha:\Theta_9\longrightarrow KO_9\cong\Z/2.
\]
To compute this index from a stable homotopy class, let
$\Omega_*^{\mathrm{fr}}$ and $\Omega_*^{\mathrm{Spin}}$ denote the bordism
rings of closed stably framed and spin manifolds, respectively.  Addition
is disjoint union and multiplication is Cartesian product; in each case
bordisms carry the same structure.  A stable normal framing induces a
spin structure on the tangent bundle.  The Pontryagin--Thom isomorphism
and the map that forgets the framing therefore define the graded-ring
homomorphism
\begin{equation}\label{eq:framed-alpha-map}
 \alpha^{\mathrm{fr}}:\pi_*^s\cong\Omega_*^{\mathrm{fr}}
 \longrightarrow\Omega_*^{\mathrm{Spin}}
 \xrightarrow{\,\alpha\,}KO_*
\end{equation}
\cite[Section 2.3, equation (2.11)]{CrowleySchick}.
The subscript on $\alpha_n^{\mathrm{fr}}$ will denote its degree-$n$ part.

In degrees eight and nine this map factors through the stable
$J$-cokernel.  Indeed, $\pi_8^s$ is finite and $KO_8\cong\Z$, so
$\alpha_8^{\mathrm{fr}}=0$.  The generator of $\im J_9$ in
\eqref{eq:pi9stable} is $\eta(\bar\nu+\epsilon)$; multiplicativity then
gives $\alpha_9^{\mathrm{fr}}(\im J_9)=0$.  Thus there are induced maps
\begin{equation}\label{eq:alpha-cokernel-factorization}
 \bar\alpha_n:\Coker J_n\longrightarrow KO_n,\qquad
 \bar\alpha_n([x])=\alpha_n^{\mathrm{fr}}(x),\qquad
 \alpha(\Sigma)=\bar\alpha_n\bigl(p_n(\Sigma)\bigr),\quad n=8,9.
\end{equation}
For the last identity, represent $p_n(\Sigma)$ by the Pontryagin--Thom
class of a stable normal framing of $\Sigma$, as in the preceding
subsection.  Forgetting that framing gives the unique spin structure
compatible with its orientation.  Changing the framing changes the
stable class by an element of $\im J_n$, whose index has just been shown
to vanish.

In particular,
$\alpha_9^{\mathrm{fr}}(\eta\epsilon)
 =\alpha_1^{\mathrm{fr}}(\eta)\alpha_8^{\mathrm{fr}}(\epsilon)=0$.
The Adams element $\mu$ satisfies
$\alpha_9^{\mathrm{fr}}(\mu)=1\in KO_9\cong\Z/2$
\cite[Section 2.3]{CrowleySchick}.  Also $\alpha(\kappa)=0$, since
$\kappa$ bounds a parallelizable manifold.  Applying
\eqref{eq:alpha-cokernel-factorization} to the basis
\eqref{eq:intro-decomp} gives
\begin{equation}\label{eq:alpha-kernel}
  \ker\alpha=\langle\kappa,\lambda\rangle
  =\{0,\kappa,\lambda,\kappa+\lambda\}.
\end{equation}

By Hitchin's obstruction, the four classes in
$\widehat\mu+\ker\alpha$ admit no metric of positive scalar curvature.
The following proposition gives the corresponding obstruction for
nonnegative sectional curvature.

\begin{proposition}[The nonnegative-sectional-curvature obstruction]
\label{prop:no-nonnegative}
If a smooth homotopy $9$-sphere $\Sigma$ has $\alpha(\Sigma)\neq0$,
then it admits no Riemannian metric of nonnegative sectional curvature.
\end{proposition}
\begin{proof}
This implication is also recorded in \cite[Proposition 7.1]{CNZ}.
Suppose $\operatorname{sec}_g\geq0$.  Let $D_g$ be the spin Dirac
operator.  An invertible Dirac operator has zero real index, so
$\alpha(\Sigma)\ne0$ implies that $\ker D_g$ contains a nonzero
harmonic spinor $\psi$.  The Lichnerowicz formula \cite{Hitchin} gives
\[
 0=\int_\Sigma\left(|\nabla\psi|^2+
                  \frac{\operatorname{Scal}_g}{4}|\psi|^2\right)dV_g.
\]
Both terms are nonnegative.  Hence $\psi$ is parallel with constant
nonzero norm, and $\operatorname{Scal}_g$ vanishes identically.
Complete an orthonormal basis of any prescribed two-plane to a tangent
orthonormal basis.  The identity
$\operatorname{Scal}_g=2\sum_{i<j}\operatorname{sec}_g(e_i\wedge e_j)$
then forces that plane's sectional curvature to vanish.  Thus $g$ is
flat.  Completeness and simple connectivity would make $(\Sigma,g)$
isometric to $\mathbb R^9$, contradicting compactness.
\end{proof}

Accordingly, the proof of \cref{thm:main-intro} reduces to constructing
positive-sectional-curvature metrics on $\kappa$, $\lambda$, and
$\kappa+\lambda$.

\subsection{Gluing discs with fixed parametrizations}
For each smooth $n$-disc $D_i$, fix a diffeomorphism
$\chi_i:D^n\to D_i$ and give $D_i$ the transported standard orientation.
The restriction $\partial\chi_i:S^{n-1}\to\partial D_i$ parametrizes
its boundary.

Given an orientation-preserving diffeomorphism
$\varphi:S^{n-1}\to S^{n-1}$, glue the two discs by identifying
$\chi_N(y)$ with $\chi_S(\varphi(y))$.  The southern disc
has the opposite orientation, as in the clutching convention above.
For $n=9$ we write
\begin{equation}\label{eq:disc-gluing}
 \Sigma=D_N^9\cup_\varphi D_S^9,
\end{equation}
with this orientation convention understood.  We also use $\varphi$
for the actual boundary map
$\partial\chi_S\circ\varphi\circ(\partial\chi_N)^{-1}$ when the
domain and target are the boundaries of the parametrized discs.

For metrics $g_i$ on $D_i$, let $h_i$ and $B_i$ be their boundary
metrics and outward second fundamental forms, expressed on the common
$S^{n-1}$ by $\partial\chi_i$.  Thus the boundary isometry condition
is exactly $h_N=\varphi^*h_S$.  Replacing $\chi_i$ by $\chi_i\circ u_i$,
where $u_i:D^n\to D^n$ is an orientation-preserving diffeomorphism,
replaces the coordinate attaching map by
\[
 \varphi'=(\partial u_S)^{-1}\circ\varphi\circ\partial u_N.
\]
Put $v_i=\partial u_i$.  The boundary tensors in the new coordinates
are $h_i'=v_i^*h_i$ and $B_i'=v_i^*B_i$.  Since
$v_S\circ\varphi'=\varphi\circ v_N$, functoriality of pullback gives
\[
 h_N'-\varphi'^*h_S'=v_N^*(h_N-\varphi^*h_S),\qquad
 B_N'+\varphi'^*B_S'=v_N^*(B_N+\varphi^*B_S).
\]
Thus changing the disc parametrizations preserves both gluing conditions.

\subsection{Group actions and relative diffeomorphisms}
For a Euclidean space $W$, put
$S(W)=\{x\in W:|x|=1\}$ and $D(W)=\{x\in W:|x|\leq1\}$.
For $r>0$, write $D^j(r)=\{x\in\R^j:|x|\leq r\}$ and
$D^j=D^j(1)$.  The symbols $\partial$ and $\Int$ denote boundary and
interior.  On Euclidean spaces, $\langle\cdot,\cdot\rangle$ denotes the
specified Euclidean inner product.  For a linear map $K$ between
inner-product spaces, adjoints
and operator norms use those inner products; in particular $K^*$ satisfies
$\langle K\xi,Y\rangle=\langle\xi,K^*Y\rangle$.
For vectors $X,Y$ we use the exterior-product norm
$|X\wedge Y|^2=|X|^2|Y|^2-\langle X,Y\rangle^2$.  We identify $S^3$ with the unit quaternions and
its Lie algebra with $\im\Hh$, with $[\xi,\zeta]=\xi\zeta-\zeta\xi$.
For $q\in S^3$, write $L_q(z)=qz$, $R_q(z)=zq$ and
$\Ad_q\xi=q\xi q^{-1}$.  The ordered basis $(i,j,k)$ fixes the Euclidean metric and orientation
of $\im\Hh$.  All metrics denoted by $h_{S(W)}$ are unit round metrics.
If $\rho:S^3\to O(W)$ is an orthogonal representation, a map
$\beta:S(W)\to S^3$ is \emph{invariant} if
$\beta(\rho(q)y)=\beta(y)$, and \emph{conjugation equivariant} if
$\beta(\rho(q)y)=q\beta(y)q^{-1}$.  The associated map
$F_\beta:S(W)\to S(W)$ is defined by $F_\beta(y)=\rho(\beta(y))y$.
We write $J_\beta$ in the equivariant case and $L_\beta$ in the invariant case.

For a compact smooth manifold $M$ with boundary, let
$\Diff_\partial(M)=\Diff(M,\mathrm{rel}\,\partial)$ be the group of
smooth diffeomorphisms equal to the identity near $\partial M$, with
composition and the $C^\infty$ topology.  A collar is a smooth embedding
$c:\partial M\times[0,\varepsilon)\to M$ with $c(x,0)=x$.
Each element of $\Diff_\partial(M)$ fixes a collar pointwise.  For each compact parameter family we
choose a common collar, fixed pointwise for diffeomorphisms and mapped
to $1$ for relative maps into $S^3$.  Relative collar compression and
smoothing give such representatives for the homotopy classes used below.
A map described as smooth on a manifold with boundary is understood
to be smooth up to the boundary.
Two relative diffeomorphisms $f_0,f_1$ of $M$ are pseudo-isotopic if
there is a diffeomorphism $F:M\times[0,1]\to M\times[0,1]$ such that
$F(x,i)=(f_i(x),i)$ for $i=0,1$ and $F$ is the identity near
$\partial M\times[0,1]$.  Unlike an isotopy, $F$ need not preserve the
intermediate slices.  We denote the group of these classes by
$\widetilde\pi_0\Diff(M,\mathrm{rel}\,\partial)$.
The basepoint of a diffeomorphism group is its identity.
For based spaces $X,Y$, $\Map_*(X,Y)$ denotes the space of continuous
basepoint-preserving maps, with the compact-open topology.  We write
$\Omega^jY=\Map_*(I^j/\partial I^j,Y)$, where $I=[0,1]$;
adjunction gives $\pi_i(\Omega^jY)=\pi_{i+j}(Y)$ at the constant map.
Other components and base maps will be specified explicitly.
The smash product is
$X\wedge Y=(X\times Y)/(X\times\{*\}\cup\{*\}\times Y)$.
For a topological group $G$, the notation $EG\to BG$ denotes a
universal principal bundle with contractible total space $EG$.

\begin{lemma}[Diffeomorphisms associated with group actions]\label{lem:associated-diffeomorphism}
Let $V$ be an orthogonal $S^3$-representation with action $\rho$.
Suppose $\beta:S(V)\to S^3$ is smooth and either invariant or
conjugation equivariant.  Then
\[
 F_\beta(y)=\rho(\beta(y))y
\]
is a smooth diffeomorphism with inverse $F_{\beta^{-1}}$, where
$\beta^{-1}(y)=\beta(y)^{-1}$.  Its signed
Jacobian for the round volume form is identically one.  A smooth homotopy through invariant or conjugation-equivariant maps
gives an isotopy of the corresponding diffeomorphisms.
\end{lemma}

\begin{proof}
In either case $\beta(F_\beta(y))=\beta(y)$; in the conjugation-equivariant
case this is the identity
$\beta(y)\beta(y)\beta(y)^{-1}=\beta(y)$.
This proves the inverse formula.
Set $K_y\xi=\left.\frac{d}{ds}\right|_{s=0}\rho(e^{s\xi})y$ and
$A_y=\beta(y)^{-1}d\beta_y$.  Differentiation gives
\[
 dF_\beta|_y=\rho(\beta(y))\bigl(I+K_yA_y\bigr).
\]
The orthogonal factor preserves the oriented round volume, since $S^3$
is connected.  By the rectangular determinant identity,
\[
 \det(I+K_yA_y)=\det(I+A_yK_y).
\]
If $\beta$ is invariant, $A_yK_y=0$.  If $\beta$ is conjugation equivariant,
\[
 A_yK_y=\Ad_{\beta(y)^{-1}}-I.
\]
The determinant is $1$ in either case, and the inverse formula is smooth
in a homotopy parameter.
\end{proof}

\subsection{Bundle conventions}
A principal right $G$-bundle $\pi:P\to B$ has a free right action
$R_h(p)=ph$ that is transitive on each fiber.  Its local trivializations
are equivariant, with action $(x,v)h=(x,vh)$.  The pullback of a bundle
$p:E\to Y$ by a map $f:X\to Y$ is
$f^*E=\{(x,e)\in X\times E:f(x)=p(e)\}$, with projection to $X$.

For $G=S^3$, a principal connection is an $\im\Hh$-valued one-form
$\omega$ on $P$ satisfying
$\omega(\xi^\#)=\xi$ and $R_h^*\omega=\Ad_{h^{-1}}\omega$, where
$\xi^\#_p=\left.\frac{d}{dt}\right|_{t=0}p\exp(t\xi)$.
Its horizontal space is $\ker\omega$, and its vertical space is
$\ker d\pi$.
For a local section $s_i$ write $A_i=s_i^*\omega$; in coordinates
$p=s_i(x)v_i$,
\[
 \omega_i=\Ad_{v_i^{-1}}A_i+v_i^{-1}dv_i.
\]
The curvature is $d\omega+\frac12[\omega\wedge\omega]$, where
$[\omega\wedge\omega](U,V)=2[\omega(U),\omega(V)]$.
It is an equivariant horizontal form and hence corresponds to a
$2$-form on $B$ with values in the adjoint bundle
\[
 \operatorname{ad}P=(P\times\im\Hh)/\bigl((p,\xi)\sim
                   (ph,\Ad_{h^{-1}}\xi)\bigr).
\]
An auxiliary left group action on a total space will be denoted by
$\star$; its quotient $P/G_\star$ is the space of its orbits
$[p]_\star=\{g\star p:g\in G\}$.

\section{Diffeomorphisms supported in a tubular neighborhood}\label{sec:tubular-neighborhood}

\subsection{Left and right quaternionic multiplication}
Fix $0<c_0<c_1<1$ and give $\R^5=\R\oplus\Hh$ the product
orientation.  Choose a smooth map
\begin{equation}\label{eq:b-definition}
 b:\R^5\longrightarrow S^3,\qquad b(u)=1\quad\text{for }|u|\geq c_0,
\end{equation}
whose restriction to $D^5(c_1)$ descends to a based map on
$D^5(c_1)/\partial D^5(c_1)\cong S^5$ representing
\begin{equation}\label{eq:b-class}
 [b]=\eta_3\eta_4\in\pi_5(S^3)\cong\Z/2.
\end{equation}
Here $\eta_2:S^3\to S^2\subset\im\Hh$ is the Hopf map
$\eta_2(q)=qiq^{-1}$, and $\eta_j=\Sigma^{j-2}\eta_2$ for $j\geq2$;
products of maps are composed from right to left.  The nonzero class in
\eqref{eq:b-class} is the classical calculation in \cite{Toda}.
A representative can be made constant on a boundary collar by a
smooth radial collapse and relative smoothing, and then extended by
$1$ to $\R^5$.  The oriented radial identification with the unit disc
fixes the identification of its relative class with $\pi_5(S^3)$.

Put $T=S^3\times D^5(c_1)$ and fix this $b$.  Smooth homotopies of $b$
that are equal to $1$ on a common boundary neighborhood induce
isotopies of the diffeomorphisms defined below.  Collar straightening identifies these relative pseudo-isotopy classes
with those used by Krylov \cite{Krylov}.
Define relative diffeomorphisms
\begin{equation}\label{eq:R-L}
  R(q,u)=(q b(u),u),
  \qquad
  L(q,u)=(b(u)q,u).
\end{equation}
They are the identity near $\partial T$ and satisfy $RL=LR$.

Let
\begin{equation}\label{eq:E-def}
  \cE=\widetilde\pi_0\Diff(T,\mathrm{rel}\,\partial)
\end{equation}
be the group of relative pseudo-isotopy classes.  Composition in $\cE$
is written multiplicatively; we also write $R,L$ for the classes of the
two diffeomorphisms in \eqref{eq:R-L}.

Embed $T$ as a tubular neighborhood in $S^8$ by
\begin{equation}\label{eq:e0}
  e_0(q,u)=\big(\sqrt{1-\abs{u}^2}\,q,u\big)
  \in S(\Hh\oplus\R^5),
\end{equation}
where $|u|\leq c_1<1$.  Its inverse recovers $q=w/|w|$ and $u$; the
quaternionic coordinate is nonzero on $e_0(T)$.  For
$f\in\Diff_\partial(T)$ define
\[
 \widehat f(x)=
 \begin{cases}
 e_0\circ f\circ e_0^{-1}(x),&x\in e_0(T),\\
 x,&x\notin e_0(T).
 \end{cases}
\]
The two formulas agree on a neighborhood of $\partial e_0(T)$.
Thus $\widehat f$ is a smooth orientation-preserving diffeomorphism.
Gluing two discs along $\widehat f$ defines the homomorphism
\begin{equation}\label{eq:Phi}
 \Phi:\cE\longrightarrow\Theta_9,\qquad
 \Phi([f])=[\Sigma_{\widehat f}].
\end{equation}
To check that $\Phi$ is well-defined, extend a relative pseudo-isotopy
to $H:S^8\times I\to S^8\times I$, with endpoint maps
$\widehat f_0,\widehat f_1$.  The diffeomorphism
$(\widehat f_1^{-1}\times\id_I)\circ H$ restricts to
$\widehat f_1^{-1}\widehat f_0$ at the bottom and to the identity at the
top.  Straighten it to the identity on a top collar and cap that end
by a $9$-disc, extending by the identity.  The capped cylinder is a
$9$-disc, so $\widehat f_1^{-1}\widehat f_0$ extends over it.  The
clutching convention of Section~\ref{sec:classification} then gives
$[\Sigma_{\widehat f_0}]=[\Sigma_{\widehat f_1}]$.
Extension preserves composition, which represents connected sum under
this convention; hence $\Phi$ is a homomorphism.

Choose an orientation-preserving embedding $e:D^8\hookrightarrow\Int T$.
The disc-supported clutching identification
$\widetilde\pi_0\Diff_\partial(D^8)\cong\Theta_9$ lets us represent
$w\in\Theta_9$ by $f_w\in\Diff_\partial(D^8)$ \cite{KM,Cerf}.
Extend $e f_w e^{-1}$ by the identity outside $e(D^8)$ and denote its
class in $\cE$ by $\iota(w)$.  Relative pseudo-isotopies extend in the
same way, so this defines a homomorphism
\begin{equation}\label{eq:iota}
 \iota:\Theta_9\longrightarrow\cE.
\end{equation}
It is injective because, with consistent orientation conventions,
\begin{equation}\label{eq:Phi-iota}
  \Phi\circ\iota=\id_{\Theta_9}.
\end{equation}
Indeed, extension first to $T$ and then to $S^8$ gives the disc-supported
clutching representative of the same oriented sphere.

\begin{lemma}[Conjugation on diffeomorphisms supported in a disc]\label{lem:A-local}
The subgroup $\cN=\iota(\Theta_9)$ is central in $\cE$.  More generally,
conjugation by any orientation-preserving diffeomorphism of $T$ fixes
$\cN$ pointwise, even when that diffeomorphism is not fixed on
$\partial T$.
\end{lemma}
\begin{proof}
Represent an element of $\cN$ by a disc diffeomorphism supported in
the image of an orientation-preserving embedding
$e:D^8\hookrightarrow\Int T$.
Conjugation by $g$ replaces $e$ by $g\circ e$ and leaves the disc
representative itself unchanged.

Shrink both parametrized embeddings to sufficiently small concentric
discs.  In coordinate neighborhoods each can be deformed to its
orientation-preserving derivative frame.  The oriented frame bundle of
the connected manifold $\Int T$ is connected, so a path of small framed
discs connects the two embeddings.  The whole path has compact image in
$\Int T$.  Isotopy extension gives an ambient isotopy fixed near
$\partial T$.  Conjugating the supported diffeomorphism along that
isotopy preserves its relative mapping class.  Thus conjugation fixes
$\cN$ pointwise.  Taking $g$ to be any relative diffeomorphism proves
centrality.
\end{proof}

\begin{lemma}\label{lem:R-L-extend}
The diffeomorphisms $\widehat R,\widehat L:S^8\to S^8$ both extend
over $D^9$.  In particular,
\[
   \Phi(R)=\Phi(L)=0.
\]
\end{lemma}

\begin{proof}
Regard
\[
 D^9=D(\Hh\oplus\R^5)
 =\{(w,u):\abs{w}^2+\abs{u}^2\leq1\}.
\]
The maps
\[
  (w,u)\longmapsto(wb(u),u),
  \qquad
  (w,u)\longmapsto(b(u)w,u)
\]
are smooth norm-preserving diffeomorphisms with explicit inverses obtained
by replacing $b$ with $b^{-1}$.  Their restrictions to $S^8$ are the
identity extensions of $R$ and $L$ under \eqref{eq:e0}.
\end{proof}

\subsection{The relative mapping-class quotient}
For $q>2$, let $C_n^q$ denote the group of smooth isotopy classes of
parametrized embeddings $S^n\hookrightarrow S^{n+q}$, with connected
sum as operation.  A normal framing of an embedding $e$ is an oriented
trivialization of its normal bundle
$e^*TS^{n+q}/TS^n$.  The group $FC_n^q$ consists of embeddings equipped
with such a framing, up to isotopy through framed embeddings; its
operation is framed connected sum.  The superscript records codimension.
Equivalently, a framed embedding is given by an oriented tubular
embedding $S^n\times D^q\hookrightarrow S^{n+q}$.  Concordance and
isotopy give the same classes in this range
\cite[Theorem 1.2 and Section 5]{Haefliger}.
Since $\cN$ is central by \cref{lem:A-local}, we may form the quotient
$\cQ=\cE/\cN$.  Write $\overline f=f\cN$ for the image of
$f\in\cE$ in $\cQ$.

\begin{theorem}[Relative mapping classes and a cardinality bound]
\label{thm:Krylov}
The quotient by $\cN=\iota(\Theta_9)$ satisfies
\begin{equation}\label{eq:Krylov}
 \cQ:=\cE/\cN\cong FC_5^4,\qquad |\cQ|\leq8.
\end{equation}
\end{theorem}
\begin{proof}
For $(p,q)=(3,5)$, Krylov's relative classification
\cite[Lemma 1 and its proof]{Krylov} gives
$\cE\cong\Theta_9\oplus FC_5^4$.  Its $\Theta_9$ summand is the
image of $\iota$, represented by diffeomorphisms supported in a disc.  Quotienting by this
summand gives $\cQ\cong FC_5^4$.

Haefliger's calculation gives $C_5^4\cong\Z/2$
\cite[Corollary 8.14]{Haefliger}.  Fix a representative of either
unframed class.  If its oriented normal bundle is not trivial, it has
no framing.  Otherwise, relative to one framing, the homotopy classes
of oriented normal framings form a torsor for
\[
 [S^5,SO(4)]=\pi_5SO(4)
 \cong\pi_5(S^3\times S^3)\cong(\Z/2)^2.
\]
Here the covering $S^3\times S^3\to SO(4)$ induces the isomorphism.
An isotopy of the underlying embedding transports its normal framing
to the chosen representative.  Isotopy extension and uniqueness of tubular
neighborhoods identify these data with tubular embeddings
$S^5\times D^4\hookrightarrow S^9$.  Self-isotopies may identify
different framings.
There are therefore at most four framed classes over each of the two
unframed classes, proving the bound.
\end{proof}

\subsection{Changing the normal framing}
Identify $\R^5=\R\oplus\Hh$ and define
\begin{equation}\label{eq:tau1}
 \tau_1:S^3\longrightarrow SO(5),\qquad
 \tau_1(q)=I_{\R}\oplus L_q.
\end{equation}
The based map $\tau_1$ changes the normal coordinates through the
diffeomorphism
\begin{equation}\label{eq:Ttau}
 T_\tau(q,u)=(q,\tau_1(q)u).
\end{equation}
Conjugation by $T_\tau$ preserves $\Diff(T,\mathrm{rel}\,\partial)$ and induces
\begin{equation}\label{eq:A-def}
 \cA:\cE\longrightarrow\cE,\qquad
 \cA[f]=[T_\tau fT_\tau^{-1}].
\end{equation}
The map $T_\tau$ preserves orientation: its derivative has diagonal
blocks $I$ and $\tau_1(q)\in SO(5)$.  Lemma~\ref{lem:A-local} therefore gives $\cA\circ\iota=\iota$.
It induces an automorphism $\overline{\cA}$ of $\cQ$, defined by
$\overline{\cA}(\overline f)=\overline{\cA(f)}$.

A second framing uses $\R^5=\R^2\oplus\im\Hh$:
\begin{equation}\label{eq:tau2}
 \tau_2(q)=I_{\R^2}\oplus\Ad_q\in SO(5).
\end{equation}
\begin{lemma}\label{lem:framing-double}
In $\pi_3SO(5)\cong\Z$,
\begin{equation}\label{eq:tau2-double}
 [\tau_2]=2[\tau_1].
\end{equation}
Consequently conjugation by $(q,u)\mapsto(q,\tau_2(q)u)$ induces
$\cA^2$ on $\cE$.
\end{lemma}
\begin{proof}
Identify $\R^5=\R\oplus\Hh=\R\oplus\R\oplus\im\Hh$ in the
stated coordinate order.  Let $c(z)=\overline z$ and put
$J_0=(-I_{\R})\oplus c\in SO(5)$.  Both factors reverse orientation,
so their direct sum preserves it.  Right multiplication by
$\overline q$ is $cL_qc$, and hence
\[
 J_0\tau_1(q)J_0^{-1}=I_{\R}\oplus R_{\overline q},\qquad
 \tau_1(q)J_0\tau_1(q)J_0^{-1}=\tau_2(q).
\]
Here $R_{\overline q}(z)=z\overline q$.
Choose $U_t\in SO(5)$, $0\leq t\leq1$, by rotation through angle
$\pi t$ in each of the ordered orthogonal planes
$\Span\{(1,0),(0,i)\}$ and $\Span\{(0,j),(0,k)\}$, fixing $(0,1)$.
Then $U_0=I$ and $U_1=J_0$.  The formula
\[
 H_t(q)=\tau_1(q)U_t\tau_1(q)U_t^{-1}
\]
defines a smooth based homotopy in $SO(5)$, with
\[
 H_0(q)=\tau_1(q)^2,\qquad H_1(q)=\tau_2(q),\qquad H_t(1)=I.
\]
Pointwise multiplication represents addition in $\pi_3SO(5)$, proving
\eqref{eq:tau2-double} with the indicated sign.

The maps $T_t(q,u)=(q,H_t(q)u)$ are diffeomorphisms of $T$ preserving
$|u|$ and fixing $S^3\times\{0\}$.  For every relative diffeomorphism $f$, the
family $T_tfT_t^{-1}$ is an isotopy relative to a common boundary collar.
At $t=0$ it is conjugation by $T_\tau^2$; at $t=1$ it is conjugation
by $(q,u)\mapsto(q,\tau_2(q)u)$.  Their actions on $\cE$ therefore agree.
\end{proof}

\subsection{Homotopy classes relative to the boundary}
For this relative-homotopy computation, identify the unit $D^5$ with
$D^5(c_1)$ by $u\mapsto c_1u$ and transport the maps by this
identification.  We retain the symbols $T,b,R,L$ in these rescaled
coordinates.  Let $\cH$ be the homotopy classes of continuous self-maps
of $T=S^3\times D^5$ that restrict to the identity on $\partial T$,
where every map in the homotopy has this same boundary restriction.
Composition initially makes $\cH$ a monoid.

\begin{lemma}[The relative-homotopy quotient]\label{lem:relative-homotopy}
The four classes of $\cH$ are represented by $\id,R,L,LR^{-1}$,
and their composition group is $(\Z/2)^2$.
Forgetting smoothness gives a surjective homomorphism
\begin{equation}\label{eq:forgetful-quotient}
 \cQ\longrightarrow\cH,\qquad
 \cK:=\ker(\cQ\to\cH),\qquad |\cK|\leq2.
\end{equation}
\end{lemma}
\begin{proof}
Write a relative map as $f=(F,G)$.  Convex interpolation of $G$ with
$u$ stays in $D^5$ and is fixed on the boundary.  Define
\begin{equation}\label{eq:relative-homotopy-coordinates}
 e_f(u)=F(1,u),\qquad
 a_f(q,u)=e_f(u)^{-1}q^{-1}F(q,u).
\end{equation}
The first component is a based map on $D^5/\partial D^5=S^5$.
The second is constant on $S^3\times\partial D^5$ and
$\{1\}\times D^5$, hence gives a based map on $S^3\wedge S^5=S^8$.
Conversely $F(q,u)=q e_f(u)a_f(q,u)$ reconstructs $F$.
These constructions apply to homotopies as well, giving a set bijection
\[
 \cH\longleftrightarrow\pi_5(S^3)\times\pi_8(S^3)
                  \cong(\Z/2)^2.
\]

Put $\beta=[b]$ and $\chi=\langle\beta,\iota_3\rangle\in\pi_8(S^3)$,
where $\iota_3$ is the class of the identity map of $S^3$.
For based maps $a:S^p\to G$ and $d:S^q\to G$ into a topological group,
the Samelson product $\langle[a],[d]\rangle\in\pi_{p+q}(G)$ is the
class induced on $S^p\wedge S^q$ by
$(x,y)\mapsto a(x)d(y)a(x)^{-1}d(y)^{-1}$.  The commutator is the
identity on the wedge, so this is a based map.  Both $\beta$ and $\chi$
are nonzero.
For the second assertion, naturality gives
\begin{equation}\label{eq:samelson-product}
 \chi=\pm\langle\iota_3,\iota_3\rangle\circ\eta_6\eta_7.
\end{equation}
The Samelson square generates $\pi_6(S^3)\cong\Z/12$
\cite[proof of Proposition III.1]{Wockel}; its $2$-primary part is an
odd multiple of Toda's order-four class $\nu'\in\pi_6(S^3)$.
The nonzero composite
\begin{equation}\label{eq:Toda-nonzero}
 \nu'\eta_6\eta_7\ne0
\end{equation}
is \cite[Proposition 5.9]{Toda}.  The pointwise coordinates of the
maps $L$ and $LR^{-1}$ are
\[
\begin{aligned}
 e_L(u)&=b(u),& a_L(q,u)&=b(u)^{-1}q^{-1}b(u)q,\\
 e_{LR^{-1}}(u)&=1,& a_{LR^{-1}}(q,u)&=q^{-1}b(u)q b(u)^{-1}.
\end{aligned}
\]
The two displayed second coordinates are commutators of based maps.
Inversion represents negation in their homotopy groups, and interchanging the smash factors
changes only the sign of the resulting class.  Bilinearity of the
Samelson product therefore identifies both second coordinates with
$\chi$, since $\pi_8(S^3)$ has order two.  Together with
$e_R=b$ and $a_R=1$, this gives
\[
\begin{array}{c|cccc}
 f&\id&R&L&LR^{-1}\\\hline
 ([e_f],[a_f])&(0,0)&(\beta,0)&(\beta,\chi)&(0,\chi).
\end{array}
\]
The maps $R,L$ commute exactly, and pointwise multiplication represents
addition in $\pi_5(S^3)$.  Since this group has order two, $b^2$ admits a
null homotopy through maps equal to $1$ on a common boundary neighborhood.  Thus $R,L$ are
involutions in $\cE$, and their four products exhaust $\cH$ and
identify its composition law with $(\Z/2)^2$.

A relative pseudo-isotopy induces a relative homotopy by projecting
its target $T\times I$ onto $T$.  If $f_w\in\Diff_\partial(D^8)$
represents $w\in\Theta_9$, the maps
\[
 h_t(x)=(1-t)f_w(x)+tx,\qquad 0\leq t\leq1,
\]
stay in the convex disc and are the identity on its boundary collar.
Transport them to the supporting disc in $T$ and extend by the identity.
This relative continuous homotopy shows that $\cE\to\cH$ kills
$\cN$ and factors
through $\cQ$.  Its surjectivity and $|\cQ|\leq8$ imply
\eqref{eq:forgetful-quotient}.
\end{proof}

\begin{lemma}[The action on relative homotopy classes]
\label{lem:framing-on-homotopy}
The automorphism $\cA$ induces the identity on $\cH$.
\end{lemma}
\begin{proof}
Let $\bar b:S^5=D^5/\partial D^5\to S^3$ be the based map induced by
$b$.  Its class is $\eta_3\eta_4$.  Write $E$ for suspension and let
$\nu_5$ generate the $2$-primary part of $\pi_8(S^5)$.  Toda's relations
\[
 \eta_3E\nu'=0,\qquad \eta_4\nu_5=(E\nu')\eta_7
\]
from \cite[Lemma 5.7 and equations (5.9)]{Toda} give
$\bar b\nu_5=\eta_3\eta_4\nu_5=0$.  The odd-primary subgroup of
$\pi_8(S^5)$ also maps to zero, since $\pi_8(S^3)=\Z/2$.  Hence
\begin{equation}\label{eq:b-postcomposition-zero}
 \bar b_*:\pi_8(S^5)\longrightarrow\pi_8(S^3)
 \quad\hbox{is zero}.
\end{equation}

The rotations of $D^5$ descend to a based family
$q\mapsto r_q$ in $\Map_*(S^5,S^5)$, where
$r_q([u])=[\tau_1(q)^{-1}u]$ and $r_1=\id$.  View these mapping spaces as fivefold loop spaces, with base maps
$\id$ and $\bar b$, respectively.
Translation by a loop is a homotopy equivalence between components,
and postcomposition preserves loop concatenation.  If $*$ denotes
concatenation, then
\[
 \Omega^5\bar b(\id*\zeta)
   =\bar b*\bigl(\Omega^5\bar b(\zeta)\bigr).
\]
Together with the standard unit homotopies, this identity identifies
the induced map on $\pi_3$ at $(\id,\bar b)$ with the homomorphism in
\eqref{eq:b-postcomposition-zero}.  Thus
$q\mapsto\bar b\circ r_q$ has a based null homotopy with base map
$\bar b$.

Set
$R_\tau=T_\tau R T_\tau^{-1}$ and
$L_\tau=T_\tau L T_\tau^{-1}$, which represent $\cA R$ and $\cA L$.
Their first components are
\[
 \pr_1R_\tau(q,u)=q b(\tau_1(q)^{-1}u),\qquad
 \pr_1L_\tau(q,u)=b(\tau_1(q)^{-1}u)q.
\]
Pulling the preceding based null homotopy back to $D^5$ changes these
components to $qb(u)$ and $b(u)q$ while fixing $\partial T$.
Convex interpolation changes each second component to $u$, also
relative to $\partial T$.  Hence $R_\tau\simeq_\partial R$ and
$L_\tau\simeq_\partial L$ as continuous maps.  Since the classes of
$R,L$ generate $\cH$, the induced action on $\cH$ is the identity.
\end{proof}

\begin{lemma}[The induced action on $\cQ$]
\label{lem:quotient-action-square}
The induced automorphism $\overline\cA$ on $\cQ$ satisfies
$\overline\cA^{\,2}=\id$.
\end{lemma}
\begin{proof}
The kernel $\cK$ has order one or two and is invariant under
$\overline\cA$, so it is fixed pointwise.  Since the action on
$\cH$ is trivial, write $\overline\cA(x)=xk_x$ for $x\in\cQ$,
with $k_x\in\cK$.  Then
\[
 \overline\cA^{\,2}(x)
 =\overline\cA(x)\overline\cA(k_x)=xk_x^2=x.
\]
\end{proof}

\section[The classes kappa and lambda]{\texorpdfstring{The classes $\kappa$ and $\lambda$}{The classes kappa and lambda}}\label{sec:classes-kappa-lambda}

We identify the sphere obtained using $\tau_1$ with the Kervaire sphere.
The framing $\tau_2$ gives a class $\lambda$ whose normal invariant is
$[\eta\epsilon]$.  Both constructions use the map $b$ from
\eqref{eq:b-definition}.

\subsection{The Kervaire class}
Consider the orthogonal representation
\begin{equation}\label{eq:rho-kappa}
 V_\kappa=\R\oplus\Hh\oplus\Hh,
 \qquad
 \rho_\kappa(q)(t,w,z)=(t,qw,qz).
\end{equation}
On the open set $w\ne0$, where the action is free, write
\[
 q=\frac{w}{\abs{w}},
 \qquad
 u=(t,q^{-1}z)\in\Int D^5(1).
\]
Restricting these coordinates to $|u|\leq c_1$ gives the tubular embedding
\begin{equation}\label{eq:e-kappa}
 e_\kappa(q,t,z)=
 \big(t,\sqrt{1-t^2-\abs{z}^2}\,q,qz\big).
\end{equation}
Let $\Pi_\kappa:\Hh\oplus(\R\oplus\Hh)\to V_\kappa$ be
$\Pi_\kappa(w,(t,v))=(t,w,v)$.  This orthogonal permutation preserves
orientation, since it moves a real line past a four-dimensional block.
With the same symbol for its restriction to the unit sphere,
\begin{equation}\label{eq:e-kappa-A}
 e_\kappa=\Pi_\kappa\circ e_0\circ T_\tau.
\end{equation}
The map $\Pi_\kappa$ extends linearly over the $9$-disc.  Conjugating
an attaching map by this coordinate change therefore preserves its
oriented homotopy-sphere class.

Define a conjugation-equivariant map on $e_\kappa(T)$ by
\begin{equation}\label{eq:beta-kappa}
  \beta_\kappa(e_\kappa(q,u))=q b(u)q^{-1},
\end{equation}
and extend it by $1$ outside $e_\kappa(T)$.  The associated
clutching diffeomorphism is
\begin{equation}\label{eq:J-beta-kappa}
  J_{\beta_\kappa}(y)=\rho_\kappa(\beta_\kappa(y))y.
\end{equation}
For $y=e_\kappa(q,u)$, equivariance of $e_\kappa$ gives
$\rho_\kappa(h)e_\kappa(q,u)=e_\kappa(hq,u)$.  Taking
$h=qb(u)q^{-1}$ therefore yields
\[
 e_\kappa^{-1}\circ J_{\beta_\kappa}\circ e_\kappa(q,u)
 =(qb(u),u)=R(q,u).
\]
Together with \eqref{eq:e-kappa-A}, this gives
\begin{equation}\label{eq:kappa-as-AR}
  [\Sigma_{J_{\beta_\kappa}}]=\Phi(\cA R).
\end{equation}

A conjugation-equivariant map on $S(V_\kappa)$ determines a path in
the gauge group of the quaternionic Hopf bundle.  We first determine
the endpoint values and then the fundamental group of this gauge group.

On $S^7=S(\Hh\oplus\Hh)$ use the free left action
$q(w,z)=(qw,qz)$.  Its quotient is $S^4=S(\R\oplus\Hh)$, with Hopf
projection
\[
 p_H(w,z)=(|w|^2-|z|^2,\,2\overline w z).
\]
We use the equivariant-function model of its gauge group,
\[
 \cG_H=\{g\in C^\infty(S^7,S^3):g(qx)=qg(x)q^{-1}\},
\]
with pointwise multiplication, identity $g\equiv1$, and the $C^\infty$
topology.  The corresponding bundle automorphism is
$x\mapsto g(x)^{-1}x$; this convention identifies pointwise multiplication
with composition of these automorphisms.  After the radial
identification of each nonpolar latitude with $S^7$, a smooth
conjugation-equivariant map on $S(\R\oplus\Hh^2)$ determines a path
of gauge transformations.  At either pole its value lies
in the center $\{1,-1\}$.

\begin{lemma}\label{lem:center-gauge}
The constant gauge transformations $1$ and $-1$ lie in distinct path
components of $\cG_H$.
\end{lemma}

\begin{proof}
Trivialize the Hopf bundle over the two hemispheres of $S^4$, with equatorial
transition $c(q)=q$.  If a path from $1$ to $-1$ in $\cG_H$ existed, it would be
represented by maps $a_N,a_S:D^4\times I\to S^3$, constant at the two ends,
with
\[
 a_N(q,s)=q a_S(q,s)q^{-1}
 \qquad(q\in S^3).
\]
Contract the southern $D^4$ to an interior point.  This deforms $a_S$,
relative to its two end faces, to a path independent of the disc coordinate.
Since $\pi_1(S^3)=0$, that path can then be deformed, still relative to
its endpoints, to $a_S(s)=\exp(\pi s i)$.  On the common side face the
corresponding homotopy of $a_N$ is prescribed by conjugation.  The homotopy
extension property for the boundary of $D^4\times I$ extends it over the
northern cylinder.  The resulting map
on the boundary of
$D^4\times I$ has side restriction
\[
 (q,s)\longmapsto q\exp(\pi s i)q^{-1}
 =\cos(\pi s)+\sin(\pi s)q i q^{-1}.
\]
After collapsing the top and bottom, this is the suspension of the Hopf
map $q\mapsto qiq^{-1}$ and represents the nonzero element
$\eta_3\in\pi_4(S^3)$.  It therefore cannot extend over $D^4\times I$.
\end{proof}

It follows that an equivariant transition has the same value at both poles.
If this value is $-1$, multiplication by the central constant $-1$ changes
$J_\beta$ by the linear map
$(t,w,z)\mapsto(t,-w,-z)$, which is orientation preserving and extends over
$D^9$.  Thus every such transition may be normalized to a based loop in
$\cG_H$ without changing its homotopy-sphere class.  A preliminary
reparameterization makes it constant on small polar caps.  Continuous
homotopies of these equivariant maps admit relative equivariant smoothing.
More precisely, let $f:M\times I\to S^3\subset\Hh$ be such a continuous
homotopy on a compact $S^3$-manifold with action $h\cdot x$, stationary near its parameter
endpoints and smooth near the invariant closed subset to be fixed.
Relative smooth approximation gives $u:M\times I\to\Hh$, equal to $f$
on a smaller neighborhood of that subset, with $\|u-f\|_\infty<1/2$.
Using normalized Haar measure, set
\[
 u_{\mathrm{av}}(x,t)=\int_{S^3}h^{-1}u(h\cdot x,t)h\,dh.
\]
The average is smooth and equivariant, preserves the prescribed values,
and satisfies $\|u_{\mathrm{av}}-f\|_\infty<1/2$.  Thus
$u_{\mathrm{av}}/|u_{\mathrm{av}}|$ is the required smooth equivariant
homotopy.  Applied relative to the polar caps, this shows that the gauge
homotopy classes determine smooth isotopy classes of the attaching maps.

\begin{lemma}\label{lem:gauge-loop}
With the identity gauge transformation as basepoint, the gauge group
satisfies
\[
 \pi_1(\cG_H)\cong\pi_5(S^3)\cong\Z/2.
\]
In the coordinates defined by $e_\kappa$, the map $b$ in
\eqref{eq:b-definition} determines the nonzero element of $\pi_1(\cG_H)$.
\end{lemma}

\begin{proof}
Fix $x_*\in S^7$ over a point in the complementary Hopf chart and set
\[
 \operatorname{ev}_{x_*}:\cG_H\to S^3,\qquad
 g\mapsto g(x_*),\qquad
 \cG_{H,*}=\operatorname{ev}_{x_*}^{-1}(1).
\]
The evaluation map is a fibration, and its fiber has weak homotopy type
$\cG_{H,*}\simeq\Omega^4S^3$ \cite[Lemmas II.7--II.8]{Wockel}.
Relative smoothing identifies the homotopy groups of the smooth and
continuous models.  The long exact sequence gives
\[
 0=\pi_2(S^3)\longrightarrow\pi_5(S^3)
 \longrightarrow\pi_1(\cG_H)\longrightarrow\pi_1(S^3)=0.
\]
A gauge transformation supported in the other, trivializing four-disc
lies in $\cG_{H,*}$ and is represented there by a map $(D^4,\partial D^4)\to(S^3,1)$.  A loop of these transformations is
a map $(I\times D^4,\partial(I\times D^4))\to(S^3,1)$ and realizes the
left-hand homomorphism from $\pi_5(S^3)$ in the displayed exact sequence.
The following coordinate change identifies this loop with $b$ in
\eqref{eq:beta-kappa}.  On a latitude of
$S(\R\oplus\Hh^2)$ the Hopf chart $w\ne0$ has affine coordinate
$\zeta=w^{-1}z\in\Hh$.  Writing $q=w/|w|$, the section and orbit
coordinates are
\begin{equation}\label{eq:Hopf-affine-chart}
 (t,q,\zeta)\longmapsto
 \left(t,\frac{\sqrt{1-t^2}}{\sqrt{1+|\zeta|^2}}q,
          \frac{\sqrt{1-t^2}}{\sqrt{1+|\zeta|^2}}q\zeta\right).
\end{equation}
Write the full normal coordinate as $u=(t,v)\in\R\oplus\Hh$.
Its quaternionic component is
\begin{equation}\label{eq:Hopf-chart-change}
 v=q^{-1}z=\frac{\sqrt{1-t^2}\,\zeta}{\sqrt{1+|\zeta|^2}},
 \qquad
 \zeta=\frac{v}{\sqrt{1-t^2-|v|^2}}.
\end{equation}
The map $(t,\zeta)\mapsto(t,v)$ is a diffeomorphism from
$(-1,1)\times\Hh$ onto $\operatorname{Int}D^5(1)$.  Its Jacobian is
\[
 \det\frac{\partial(t,v)}{\partial(t,\zeta)}
 =\frac{(1-t^2)^2}{(1+|\zeta|^2)^3}>0.
\]
The coordinate change is an orientation-preserving diffeomorphism, so
it preserves compact support and extends to degree $+1$ on the one-point
compactifications.  The pullback of $b$ is therefore the nonzero supported
loop of gauge transformations.  The corresponding equivariant map is
$q b(t,q^{-1}z)q^{-1}$, as in \eqref{eq:beta-kappa}.
\end{proof}

\begin{proposition}\label{prop:kervaire-identification}
The sphere in \eqref{eq:kappa-as-AR} is the Kervaire sphere:
\begin{equation}\label{eq:kervaire-class}
  \Phi(\cA R)=\kappa.
\end{equation}
\end{proposition}

\begin{proof}
Let $K^9$ denote the boundary of the plumbing of two tangent disc
bundles of $S^5$, oriented so that $[K^9]=\kappa$.  This is the nonzero
class of $bP_{10}\cong\Z/2$ and hence is exotic
\cite{KM,Wang}.  We identify the
clutched sphere with this manifold using the real pullback construction
of \cite[Section 4]{Speranca}.  Write
$\R^6=\R e_0\oplus\R^5\cong\mathbb C^3$, with $ie_0=e_1\in\R^5$.
On $S^9\subset\R^5\oplus\R^5$ consider
\[
 f(x,y)=\exp_{e_0}\!\left(\pi\tau(y/|y|)x\right),
 \qquad \tau(v)w=2\langle v,w\rangle v-w.
\]
Here $\exp_{e_0}$ is the exponential map of the unit round $S^5$,
and $\tau(v)$ is reflection in the line spanned by $v$.  Although the
formula uses $y/|y|$, the map extends smoothly across $y=0$.
With $a=|x|^2$, it can be written as
\begin{align*}
 f(x,y)={}&\cos(\pi\sqrt a)e_0
 -\frac{\sin(\pi\sqrt a)}{\sqrt a}x\\*
 &+\frac{2\sin(\pi\sqrt a)}{\sqrt a(1-a)}\langle x,y\rangle y.
\end{align*}
Each scalar coefficient extends smoothly to $0\leq a\leq1$: at
$a=0$ use the power series of $\sin(\pi\sqrt a)/\sqrt a$, and at
$a=1$ use its simple zero.  Thus $f$ is smooth.  Orthogonal transformations commute with both
the reflection formula and the round exponential map, so $f$ is
$SO(5)$-equivariant.  It has norm one where the exponential formula is
defined, and hence everywhere by continuity.

Use the first-column projections $SO(6)\to S^5$, $A\mapsto Ae_0$,
and $SU(3)\to S^5$ given by the same formula.  Thus
$E=f^*SO(6)=\{(x,y,A):Ae_0=f(x,y)\}$.  The free $\star$-action is
$h\star(x,y,A)=(hx,hy,hA)$, with $SO(5)$ fixing $e_0$.
Speran\c ca's plumbing identification gives $E/SO(5)_\star\cong K^9$.
Passing from the orthogonal construction to its determinant-positive
component does not change this quotient: every $O(5)$ orbit meets that
component in precisely an $SO(5)$ orbit.

Set $P=f^*SU(3)\subset E$, where $SU(3)$ is viewed as a subgroup of $SO(6)$
through its real action on $\mathbb C^3$.  The subgroup $SO(5)$ fixes
$e_0$, and $SU(2)$ is embedded in the last two complex coordinates.
Transitivity of $SU(3)$ on $S^5$, with stabilizer $SU(2)$ at $e_0$,
gives
\[
 SO(6)=SO(5)SU(3),\qquad SO(5)\cap SU(3)=SU(2).
\]  They give a diffeomorphism
$SO(5)\times_{SU(2)}P\longrightarrow E$,
$[h,p]\mapsto h\star p$.  Here the associated space is the quotient of
$SO(5)\times P$ by $(h,p)\sim(hk^{-1},k\star p)$ for $k\in SU(2)$,
so the displayed map is well-defined.  Consequently
\begin{equation}\label{eq:SU2-reduction}
 P/SU(2)_\star\cong E/SO(5)_\star\cong K^9.
\end{equation}
The bundle $P\to S^9$ is a principal right $SU(2)$-bundle: the right
action multiplies the last two columns, while leaving the first column
fixed.  It commutes with left multiplication in the restricted $\star$-action.

On $\R^5=\R e_1\oplus\mathbb C^2$, the restricted $SU(2)$ action is
trivial on $e_1$ and is the defining complex representation on
$\mathbb C^2$.  To identify it with left quaternionic multiplication,
write $q=a+jb$, with $a,b\in\mathbb C$, and use
\[
 A(q)=\begin{pmatrix}a&-\overline b\\ b&\overline a\end{pmatrix},
 \qquad I(z_1,z_2)=z_1+jz_2.
\]
Then $I(A(q)z)=qI(z)$.  The real map $I$ sends the complex-oriented basis
to $(1,i,j,-k)$; its determinant is $-1$, so $I\oplus I$ preserves
orientation.  Thus the base representation is
$\R^2\oplus\Hh^2$, with the left action used in
\eqref{eq:rho-kappa}.  Choose poles $(e_1,0)$ and
$(-e_1,0)$, both mapped by $f$ to $-e_0$.  The same matrix
\[
 A_0=\operatorname{diag}(-1,i,i)\in SU(3)
\]
specifies a point in each polar fiber and commutes with the embedded $SU(2)$.  At either pole,
$h\star A_0=A_0\cdot h$.  Choose a principal connection and average it under the $\star$-action.
The averaged form is still a principal connection, because the two
actions commute.  Parallel transport from these two fiber points along
radial geodesics gives equivariant sections $s_i$ satisfying
$h\star s_i(y)=s_i(\rho_\kappa(h)y)h$.
Consequently the $\star$-action in both trivializations is
$(y,v)\mapsto(\rho_\kappa(h)y,hv)$.  If
$s_N(y)=s_S(y)\gamma(y)$, these identities imply
$\gamma(\rho_\kappa(h)y)=h\gamma(y)h^{-1}$.
Thus the transition between the two trivializations is conjugation equivariant.

By \cref{lem:center-gauge,lem:gauge-loop}, such transitions have at most
two homotopy-sphere classes.  The constant loop gives the standard sphere, whereas
\eqref{eq:SU2-reduction} gives the nonstandard Kervaire sphere.  The latter
therefore corresponds to the nonzero element of $\pi_1(\cG_H)$.  The map $b$
represents that element as well, proving \eqref{eq:kervaire-class}.
\end{proof}

\subsection[Construction of lambda]{\texorpdfstring{Construction of $\lambda$}{Construction of lambda}}
Set
\begin{equation}\label{eq:Vm}
 V_m=\R^m\oplus\Hh\oplus\im\Hh,
 \qquad
 \rho_m(q)(a,w,z)=(a,qw,qzq^{-1}),
 \qquad m=1,2.
\end{equation}
We use $m=1$ to identify the exotic $8$-sphere and $m=2$ to construct
$\lambda\in\Theta_9$.  On $w\ne0$, where the action is free, put
\[
 q=\frac w{\abs w},\qquad v=q^{-1}zq.
\]
Define the coordinates $\psi_m$ by
\begin{equation}\label{eq:Vm-coordinates}
 \psi_m(q,a,v)=
 \left(a,q\sqrt{1-\abs a^2-\abs v^2},qvq^{-1}\right),
  \qquad (a,v)\in\Int D^{m+3}(1).
\end{equation}
Let $c:D^{m+3}(1)\to S^3$ be smooth and equal to $1$ on a boundary
neighborhood, as in the support convention of \cref{sec:tubular-neighborhood}.
Define
\begin{equation}\label{eq:beta-c-general}
 \beta_c(a,w,z)=q c(a,q^{-1}zq)q^{-1}
\end{equation}
on $w\ne0$ and set it equal to $1$ near $w=0$.  Since $c$ is
constant on a boundary collar, the two formulas agree on an invariant
neighborhood; thus $\beta_c$ is globally smooth and conjugation equivariant.
In these coordinates, $J_{\beta_c}$ acts by right multiplication on the
$S^3$ factor:
\begin{equation}\label{eq:beta-c-coordinates}
 \psi_m^{-1}\circ J_{\beta_c}\circ\psi_m(q,a,v)=(qc(a,v),a,v).
\end{equation}

\begin{lemma}[Equivariant normalization]\label{lem:equivariant-normalization}
For $m=1,2$, every smooth conjugation-equivariant map
$\beta:S(V_m)\to S^3$, after multiplication by a possible central constant
$-1$, is equivariantly homotopic through smooth maps to some $\beta_c$ of
\eqref{eq:beta-c-general}.  Neither operation changes the oriented
homotopy-sphere class defined by $J_\beta$.
\end{lemma}

\begin{proof}
The set of points with nontrivial stabilizer is
$F=\{w=0\}=S(\R^m\oplus\im\Hh)$.  At a point $(a,0,z)$ with $z\ne0$ the
stabilizer is the circle fixing the axis $z$, so equivariance forces
\[
 \beta(a,0,z)=\cos\vartheta(a)
 +\sin\vartheta(a)\frac z{\abs z}
\]
for a continuous circle-valued phase on $D^m$.  A real lift exists because
$D^m$ is simply connected.  At $z=0$ the value is central.  We first contract the phase continuously and then smooth the resulting
equivariant homotopy relative to its prescribed data.  For $m=2$ the
boundary circle is connected; after multiplication by $-1$ if necessary,
the boundary value is $1$, and the phase lifts to a real-valued function
vanishing on the boundary.  For $m=1$, use the continuous section
\[
 (a,v)\longmapsto\bigl(a,\sqrt{1-a^2-|v|^2},v\bigr),
 \qquad (a,v)\in D^4.
\]
Composing with $\beta$ extends $\beta|_F$ continuously over $D^4$.
Replace its real phase by a homotopic smooth function with the same
endpoint values in $\pi\Z$, constant near those endpoints.  The round
volume form of $S^3$ in these coordinates is
$\sin^2\vartheta\,d\vartheta\wedge dV_{S^2}$.
Since $\operatorname{vol}(S^2)=4\pi$ and
$\operatorname{vol}(S^3)=2\pi^2$, the degree, up to the fixed boundary
orientation sign, is
\[
 \frac2\pi\int_{-1}^{1}\sin^2\vartheta(a)\vartheta'(a)\,\dd a
 =\frac{\vartheta(1)-\vartheta(-1)}\pi.
\]
The last equality follows from
$\int\sin^2 t\,dt=t/2-\sin(2t)/4$ and the endpoint values.
The extension over $D^4$ forces this degree to vanish, so the endpoint
values of the lift agree.  After a central
normalization and subtracting a constant multiple of $2\pi$, the phase
vanishes at both endpoints.

Contracting the phase makes $\beta|_F$ equal to $1$.  The invariant
submanifold $F$ has an equivariant tubular neighborhood, so the pair
$(S(V_m),F)$ has the equivariant homotopy extension property.  Extend the
contraction continuously over the sphere.  At its final time the map is
uniformly close to $1$ on an invariant neighborhood of $F$.  In that
neighborhood use the equivariant logarithm at $1$ and an invariant cutoff
to make the final map exactly $1$ on a smaller neighborhood.  Arrange the
homotopy to be stationary near its two time endpoints.

Apply the relative equivariant smoothing procedure following
\cref{lem:center-gauge}, fixing the original smooth endpoint and the
final constant neighborhood.  This gives a smooth equivariant homotopy
ending at a map equal to $1$ near $F$.  Its restriction to the free section is a
smooth map $c$ equal to $1$ near $\partial D^{m+3}$.

By \cref{lem:associated-diffeomorphism}, an equivariant homotopy gives an isotopy
of the maps $J_\beta$.
Multiplication by $-1$ composes $J_\beta$ with a fixed orientation-preserving
linear map that extends over the disc.  Hence neither change alters the
oriented sphere class.
\end{proof}

\subsubsection*{A homomorphism to $\Diff_\partial(D^6)$}
We construct a homomorphism that relates the maps $J_{\beta_c}$ to the
Gromoll construction.
Identify $D^3$ with the unit disc in $\im\Hh$.  For $0<\sigma<1$, let
\[
 \cG_{\mathrm c}=\{c\in C^\infty(D^3,S^3):c(v)=1\text{ for }|v|\geq\sigma\},
 \qquad \cD=\Diff_\partial(D^6),
\]
with pointwise multiplication on $\cG_{\mathrm c}$ and the $C^\infty$
topologies.  Relative smoothing and radial compression
show that inclusion in $\Map_*(D^3/\partial D^3,S^3)$ is a weak
equivalence of topological groups.  In particular,
$\pi_r\cG_{\mathrm c}\cong\pi_{r+3}(S^3)$ and
$\pi_0\cG_{\mathrm c}\cong\Z$.

Use the representation $\rho_0(q)(w,z)=(qw,qzq^{-1})$ on
$S(\Hh\oplus\im\Hh)=S^6$.  The region $w\ne0$ has coordinates
\[
 \psi_0(q,v)=\bigl(\sqrt{1-|v|^2}\,q,qvq^{-1}\bigr),
 \qquad q\in S^3,\quad v\in\Int D^3.
\]
The inverse recovers $q=w/|w|$ and $v=q^{-1}zq$; in particular
$|v|=|z|$.  Define the sphere diffeomorphism $\cT(c)$ by
\begin{equation}\label{eq:T-homomorphism}
 \cT(c)\bigl(\psi_0(q,v)\bigr)
   =\psi_0(qc(v)^{-1},v),
\end{equation}
and extend it by the identity near $w=0$.  This extension is smooth
because $c=1$ for $|v|\geq\sigma$.  All these maps fix the same open set
$\{|z|>\sigma\}$.  Choose an open round cap $C$ whose closure lies in
that set, and fix an orientation-preserving parametrization
$j:D^6\to S^6\setminus C$.  Every $\cT(c)$ is the identity on a
neighborhood of the closed cap, hence preserves its complement and is
the identity near its boundary.  Restriction and conjugation by $j$
therefore define $\cT:\cG_{\mathrm c}\to\cD$; we use the same symbol
for the sphere map and this disc map.

For $r=1,2$ define
\begin{equation}\label{eq:Lambda-r}
 \Lambda_r:\pi_r\Diff_\partial(D^6)\longrightarrow\Theta_{r+7}
\end{equation}
as follows.  Choose a smooth representative family $f_a$,
$a\in D^r$, such that $f_a(z)=z$ near the entire boundary of
$D^r\times D^6$.  Form $\widehat f(a,z)=(a,f_a(z))$ and round
the corners within its identity collar.  Extend the resulting
$D^{r+6}$-diffeomorphism by the identity to $S^{r+6}$, and clutch two
$(r+7)$-discs.  This is the Gromoll homomorphism with the clutching
convention fixed in Section~\ref{sec:classification}
\cite[Section 2.1]{CrowleySchick}.

\begin{lemma}[The Gromoll construction for $\cT$]\label{lem:T-Gromoll}
The map $\cT$ is a continuous group homomorphism.  Let $r=1$ or $2$,
and let $c:D^r\times D^3\to S^3$ be a smooth family representing a class
in $\pi_r\cG_{\mathrm c}$, equal to $1$ near the product boundary.
After compressing its support into the unit normal disc as below, the
associated attaching map satisfies
\[
 [\Sigma_{J_{\beta_c}}]=-\Lambda_r(\cT_*[c])
                         =\Lambda_r(\cT_*[c])\in\Theta_{r+7}.
\]
Here $\Theta_8$ and $\Theta_9$ have exponent two.
\end{lemma}

\begin{proof}
In the fixed coordinates $\psi_0$, applying $\cT(d)$ and then
$\cT(c)$ sends
\[
 (q,v)\longmapsto(qd(v)^{-1}c(v)^{-1},v)
              =(q(c(v)d(v))^{-1},v).
\]
Thus $\cT(c)\circ\cT(d)=\cT(cd)$, both on $S^6$ and on the fixed
disc.  The inverse is $\cT(c^{-1})$.  The common identity neighborhood
ensures smooth extension, and the coordinate formula gives continuity
in the $C^\infty$ topology.

For the family assertion, use orientation-preserving relative coordinate
compressions to arrange
\[
 \supp(c-1)\subset\{|a|<a_*,\ |v|<d_*\},\qquad
 0<a_*<1,\quad 0<d_*<\sigma\sqrt{1-a_*^2}.
\]
These compressions are isotopic to the identity and preserve the
relative class.  Extend $c$ by $1$ outside the product disc.  The support
condition implies $a_*^2+d_*^2<1$, so restricting this extension to the
unit normal disc gives the map $c$ used to define $J_{\beta_c}$.
Write a point of $S^{r+6}$ as $(a,\sqrt{1-|a|^2}\,x)$, with
$x\in S^6$ and $|a|<1$.  In the coordinates $\psi_0$, its induced
family is
\[
 x\longmapsto\cT(c_a)^{-1}(x),\qquad
 c_a(v)=c(a,\sqrt{1-|a|^2}\,v).
\]
For $|a|\leq a_*$ the normal support of $c_a$ lies in the fixed disc of
radius $d_*/\sqrt{1-a_*^2}<\sigma$; for $|a|\geq a_*$ the map is
constant.  Thus every $c_a$ lies in the same group $\cG_{\mathrm c}$ and every
$\cT(c_a)$ fixes the chosen cap.  The rescaling is a relative homotopy:
replace $\sqrt{1-|a|^2}$ by
$1-t+t\sqrt{1-|a|^2}$, $0\leq t\leq1$.  The same support estimate
holds throughout.  The family is constant near $|a|=1$ and near the
product boundary, so it remains smooth after rounding the corners.

Choose
\begin{equation}\label{eq:parameter-interior-radius}
 a_*<a_1<1.
\end{equation}
Use the fixed parametrization $j:D^6\to S^6\setminus C$ chosen
above.  The map
\begin{equation}\label{eq:parameter-product-embedding}
 E:D^r(a_1)\times D^6\longrightarrow S^{r+6},\qquad
 E(a,x)=\bigl(a,\sqrt{1-|a|^2}\,j(x)\bigr)
\end{equation}
is a smooth embedding, including its corners: the scale is at least
$\sqrt{1-a_1^2}>0$, and the inverse on its image recovers $a$ and then
$j(x)$ by division by that scale.  

The entire nontrivial support lies in the interior of this smaller
product.  On its parameter boundary collar, $|a|>a_*$ and the family is
constant; near the $D^6$ boundary, every $\cT(c_a)$ fixes the common cap.
Thus in the chart \eqref{eq:parameter-product-embedding} the sphere map is
exactly
\[
 (a,x)\longmapsto(a,\cT(c_a)^{-1}(x)),
\]
with a common identity collar on the whole product boundary.  Round the
corners within that collar, and rescale the parameter disc with the
orientation transported by the rescaling.  This gives the defining
Gromoll representative \cite[Section 1.1]{Wang}, equivalently the
cubical model with fixed boundary collars in \cite[Section 2.1]{CrowleySchick}.
Extension by the identity outside the chart recovers the original
sphere map, and inversion accounts for the minus sign in the statement.
\end{proof}

\subsubsection*{Precomposition with the Hopf map}
For a topological group $H$, precomposition with the Hopf map
$\eta_2:S^3\to S^2$ defines the following operation on pointed homotopy sets:
\begin{equation}\label{eq:Hopf-operation}
 \mathfrak h_H:\pi_1H\cong\pi_2BH
 \xrightarrow{\,\circ\eta_2\,}\pi_3BH\cong\pi_2H.
\end{equation}
It is natural under continuous group homomorphisms.  The universal principal bundle gives the comparison
$B\cG_{\mathrm c}\simeq\Map_*(S^3,BS^3)$.  Applying
$\Map_*(S^3,-)$ to the universal principal bundle
$S^3\to ES^3\to BS^3$ gives a fibration with contractible total space
and fiber $\Map_*(S^3,S^3)$.  Its base is connected since
$\pi_3(BS^3)=\pi_2(S^3)=0$.  Relative smoothing and the fixed-collar
compression identify this fiber weakly with $\cG_{\mathrm c}$.  Under
adjunction, precomposition by $\eta_2$ becomes precomposition by
$\Sigma^3\eta_2=\eta_5$ in $BS^3$.  To pass to the fiber, write
$\partial_j:\pi_j(BS^3)\to\pi_{j-1}(S^3)$ for the connecting
isomorphism.  Since $\eta_5=\Sigma\eta_4$, loop adjunction and naturality
of the map $\Omega BS^3\to S^3$ give
\[
 \partial_6(v\circ\eta_5)=(\partial_5v)\circ\eta_4,
 \qquad v\in\pi_5(BS^3).
\]
Taking $\partial_5v=\eta_3$ therefore gives
\begin{equation}\label{eq:Hopf-on-gauge}
 \mathfrak h_{\cG_{\mathrm c}}(\eta_3)=\eta_3\eta_4.
\end{equation}

\begin{lemma}[Gromoll--Morlet compatibility]\label{lem:Hopf-detection}
For every $x\in\pi_1\cD$,
\begin{equation}\label{eq:Hopf-detection}
 p_9\Lambda_2(\mathfrak h_{\cD}x)
 =\bigl(p_8\Lambda_1(x)\bigr)\eta
 \quad\text{in }\Coker J_9.
\end{equation}
\end{lemma}

\begin{proof}
Let $O_j=O(j)$.  We use the simplicial groups $PL_j$ and $TOP_j$ of
piecewise-linear and topological automorphisms of $\R^j$ fixing the
origin, in the smoothing models of \cite{SalvatoreTurchin}.
The notation $PL_j/O_j$ denotes the homotopy fiber of
$BO_j\to BPL_j$.  For a based map $f:X\to Y$, its homotopy fiber is
modeled by pairs $(x,\gamma)$ with a path $\gamma$ from $f(x)$ to the
basepoint of $Y$.  This homotopy fiber is represented by
$PD_j/O_j$, with $PD_j$ the piecewise-differentiable intermediary in
\cite[Sections 5.1, 8.1 and 9.2]{SalvatoreTurchin}.  The map
$PD_j\to TOP_j$ supplies its natural comparison with $TOP_j/O_j$.
Write $PL/O$ for the stable smoothing space and put
$X_6=PL_6/O_6$, $Y_6=TOP_6/O_6$, with these conventions throughout.
Use the smoothing comparison
$M_6:\pi_i\cD\to\pi_{i+7}X_6$ obtained from
\cite[Section 9.2]{SalvatoreTurchin}.  For a smooth relative
diffeomorphism $f$ of $D^6$, both the PD and TOP framed models assign
the derivative $Df$ and the path of germs
\[
 H^f_{z,t}(v)=
 \begin{cases}
 (f(z+tv)-f(z))/t,&t>0,\\
 Df_z(v),&t=0.
 \end{cases}
\]
The forgetful map from PD to TOP preserves these data and the fibers
over the identity.  It induces the natural map between the two
comparison diagrams in \cite[Section 9.2, after equation (9.4)]{SalvatoreTurchin};
the displayed path is equation (8.4) there.
For relative diffeomorphisms $f,g$, the composition identity is
\[
 H^{f\circ g}_{z,t}=H^f_{g(z),t}\circ H^g_{z,t}.
\]
It holds for $t>0$ by substitution and at $t=0$ by the chain rule.
Thus the TOP map preserves composition, as in Lemma 8.1.
The fiber inclusion and
Proposition A.1 then identify this comparison as an $E_1$ equivalence
$\cD\simeq\Omega^7Y_6$, where $E_1$ denotes multiplication associative
up to coherent homotopy; see the proof of Theorem 9.1 in that paper.
Delooping gives $B\cD\simeq\Omega^6_0Y_6$, with subscript $0$
denoting the component of the constant map.  By adjunction,
precomposition by $\eta_2$ corresponds to precomposition by
$\Sigma^6\eta_2=\eta_8$ in $Y_6$.

The natural comparison lies in the fiber sequence
\[
 X_6\xrightarrow{\,j\,}Y_6\longrightarrow TOP_6/PL_6
        \simeq K(\Z/2,3)
\]
from \cite[equation (2.3) and Section 9.2]{SalvatoreTurchin}.
Here $K(\Z/2,3)$ denotes the Eilenberg--MacLane space whose only nonzero
homotopy group is $\pi_3=\Z/2$.  In particular its homotopy groups in
degrees $8$, $9$, and $10$ vanish.
The long exact sequence therefore makes $j_*$ an isomorphism on
$\pi_8$ and $\pi_9$.
The preceding map between comparison diagrams identifies $j_*M_6$
with the TOP comparison.  It follows that
$j_*M_6(\mathfrak h_{\cD}x)=(j_*M_6x)\circ\eta_8$.
Postcomposition with $j$ commutes with precomposition by $\eta_8$.
Since $j_*$ is injective on $\pi_9$, it follows that
\begin{equation}\label{eq:actual-Morlet-Hopf}
 M_6(\mathfrak h_{\cD}x)=(M_6x)\circ\eta_8.
\end{equation}

For the geometric Gromoll map of \eqref{eq:Lambda-r},
\cref{lem:parameterized-comparison} identifies this smoothing
comparison with stabilization.  Let
$\Psi_j:\Theta_j\to\pi_j(PL/O)$ be the smoothing-class isomorphism and
$S:X_6\to PL/O$ the stabilization map.  Then
\begin{equation}\label{eq:actual-sphere-stabilization}
 \Psi_{r+7}\Lambda_r=S_*M_6,\qquad r=1,2.
\end{equation}
The smoothing class of the glued sphere is identified in
\cite[Lemmas 2.5 and 2.7]{CrowleySchick}.

It remains to identify the normal invariant of the clutched sphere.
Let $G$ be the stable monoid of based homotopy self-equivalences of
spheres, based at the identity.  The canonical maps $BO\to BPL\to BG$
forget a stable vector-bundle structure first to a piecewise-linear
bundle and then to its underlying stable spherical fibration.
As with the smoothing spaces above, we use homotopy-fiber models:
\[
 PL/O=\operatorname{hofiber}(BO\to BPL),\qquad
 G/O=\operatorname{hofiber}(BO\to BG).
\]
The identity on $BO$ and the map $BPL\to BG$ induce the
normal-invariant map $\mathfrak n:PL/O\to G/O$.  Under the canonical
identification $\pi_j(G)\cong\pi_j^s$ for $j>0$, the map $O\to G$
induces $J_j$.  The homotopy-fiber sequence $O\to G\to G/O$ therefore
gives, for $j=8,9$, an injection
\[
 \jmath_j:\Coker J_j\hookrightarrow\pi_j(G/O)
\]
whose image is the torsion subgroup of $\pi_j(G/O)$
\cite[Section 2.2]{CrowleySchick}.  The group
$\pi_j(PL/O)\cong\Theta_j$ is finite, so the image of
$\mathfrak n_*$ lies in this subgroup.  We may therefore define
\[
 P_j=\jmath_j^{-1}\circ\mathfrak n_*:
       \pi_j(PL/O)\longrightarrow\Coker J_j,
\]
using the inverse of $\jmath_j$ on its image.
By \cite[Proposition 2.10]{CrowleySchick},
$P_j\Psi_j=p_j$.

Precomposition by $\eta_8:S^9\to S^8$ commutes with $S$ and $\mathfrak n$.
Its descent to the two $J$-cokernels can be checked directly:
\begin{equation}\label{eq:J-descent-eight-nine}
 \im J_8=\langle\bar\nu+\epsilon\rangle,
 \qquad
 (\bar\nu+\epsilon)\eta
   =\eta\bar\nu+\eta\epsilon\in\im J_9.
\end{equation}
The image descriptions are given in \cite[Appendix C, Table 4]{Wang}.
Since $|\bar\nu|=|\epsilon|=8$, graded commutativity has sign $+1$.
Thus multiplication by $\eta$ descends to
$\Coker J_8\to\Coker J_9$.  This is the required special
case of \cite[Proposition 2.4]{Wang}, with $(n,i,t)=(6,1,1)$.
Consequently
\[
\begin{tikzcd}[column sep=large]
\pi_1\cD \arrow[r,"\mathfrak h_{\cD}"]
  \arrow[d,"M_6"'] & \pi_2\cD \arrow[d,"M_6"]\\
\pi_8X_6 \arrow[r,"\circ\eta_8"]
  \arrow[d,"P_8S_*"'] & \pi_9X_6 \arrow[d,"P_9S_*"]\\
\Coker J_8 \arrow[r,"\cdot\eta"] & \Coker J_9
\end{tikzcd}
\]
commutes.  By \eqref{eq:actual-sphere-stabilization}, its vertical
composites are $p_8\Lambda_1$ and $p_9\Lambda_2$.  This proves the lemma.
\end{proof}

\begin{lemma}[Identification of the exotic $8$-sphere]\label{lem:exotic-eight-sphere}
For the nonzero class $\eta_3\in\pi_4(S^3)=\pi_1\cG_{\mathrm c}$,
\begin{equation}\label{eq:exotic-eight-sphere}
 \Lambda_1(\cT_*\eta_3)=\Sigma^8_{[\epsilon]},
 \qquad p_8(\Sigma^8_{[\epsilon]})=[\epsilon].
\end{equation}
\end{lemma}
\begin{proof}
Speran\c ca's pullback is defined by
\[
 f_8(\ell,w,z)=
 \frac{(\ell+wi\overline w,\,z)}
      {\sqrt{\ell^2+|w|^4+|z|^2}}
 :S(\R\oplus\Hh\oplus\Hh)\longrightarrow S^7,
\]
using the first-column projection $Sp(2)\to S^7$.  On the domain,
\[
 \ell^2+|w|^4+|z|^2=1-|w|^2+|w|^4\geq\tfrac34.
\]
Since $wi\overline w$ is imaginary, the squared norm of the numerator
is exactly the same expression.  Hence $f_8$ is a globally smooth
sphere-valued map.  Let $P_8=f_8^*Sp(2)$ be the pullback bundle.
Write a point as $(x,(v_1,v_2))$, where $x=(\ell,w,z)$,
$v_1=f_8(x)$, and $(v_1,v_2)\in Sp(2)$.  Its principal right action
is $(v_1,v_2)\cdot h=(v_1,v_2h)$.  The commuting $\star$-action is
\[
 q\star\bigl((\ell,w,z),(v_1,v_2)\bigr)
 =\bigl((\ell,qw,qzq^{-1}),(qv_1q^{-1},qv_2)\bigr).
\]
Quaternionic multiplication acts componentwise on the columns and
preserves their orthonormality.  Equivariance of $f_8$ therefore makes
this an action on $P_8$.  It is free:
$qv_2=v_2$ for a nonzero column forces $q=1$.  Speran\c ca identifies
$P_8/S^3_\star$ with the unique exotic $8$-sphere
\cite[Section 2.1 and Theorem 1]{Speranca}.
At the two poles $(\pm1,0,0)$, the first columns are $(\pm1,0)$.
Choose $(0,1)$ for the second column in each fiber.  At these chosen fiber
points, the action gives $(0,q)=(0,1)q$, agreeing with
principal right multiplication by $q$.  Equivariant radial parallel transport therefore
gives two trivializations whose transition on the equator is conjugation equivariant.
Writing $z=a+z_0$ with $a=\operatorname{Re}z$ and
$z_0\in\im\Hh$ identifies that equatorial representation with
$V_1=\R\oplus\Hh\oplus\im\Hh$.  By \cref{lem:equivariant-normalization}, its transition is represented by
one of the two elements of $\pi_4(S^3)$.  The zero element has a
null homotopy equal to $1$ on a common boundary neighborhood, which gives an isotopy of the attaching map to
the identity and hence the standard sphere.  Since this quotient
is exotic, the resulting map $c$ represents $\eta_3$.
The conclusion follows from \cref{lem:T-Gromoll}.

\end{proof}

For $m=2$, put
\begin{equation}\label{eq:rho-lambda}
 V_\lambda=V_2=\R^2\oplus\Hh\oplus\im\Hh,
 \qquad
 \rho_\lambda(q)(a,w,z)=(a,qw,qzq^{-1}),
\end{equation}
and, for $q\in S^3$ and $(a,z)\in D^5(c_1)\subset\R^2\oplus\im\Hh$,
define
\begin{equation}\label{eq:e-lambda}
 e_\lambda(q,a,z)=
 \big(a,\sqrt{1-\abs a^2-\abs z^2}\,q,qzq^{-1}\big).
\end{equation}
Use the orthogonal permutation
$\Pi_\lambda:\Hh\oplus(\R^2\oplus\im\Hh)\to V_\lambda$ given by
$\Pi_\lambda(w,(a,v))=(a,w,v)$, and put
$T_{\tau_2}(q,u)=(q,\tau_2(q)u)$.  The permutation preserves orientation
and extends linearly over the $9$-disc.  The explicit normal-framing
homotopy of \cref{lem:framing-double} gives
\begin{equation}\label{eq:e-lambda-A2}
 e_\lambda=\Pi_\lambda\circ e_0\circ T_{\tau_2}
 \simeq\Pi_\lambda\circ e_0\circ T_\tau^2
\end{equation}
through parametrized tubular embeddings with fixed core.  Conjugation
by $\Pi_\lambda$ and isotopy extension preserve the clutched sphere
class.
Taking $c=b$ in
\eqref{eq:beta-c-general}, set
\begin{equation}\label{eq:beta-lambda}
 \beta_\lambda(e_\lambda(q,u))=q b(u)q^{-1}.
\end{equation}
In the coordinates $e_\lambda$, the attaching diffeomorphism is $R$, so
\begin{equation}\label{eq:lambda-A2R}
 [\Sigma_{J_{\beta_\lambda}}]=\Phi(\cA^2R).
\end{equation}
Define
\begin{equation}\label{eq:lambda-def}
 \lambda:=\Phi(\cA^2R).
\end{equation}

\begin{proposition}\label{prop:lambda-invariants}
The class $\lambda$ satisfies
\[
 p_9(\lambda)=[\eta\epsilon]=[\nu^3]\ne0,
 \qquad
 \alpha(\lambda)=0.
\]
In particular, $\lambda\notin bP_{10}$.
\end{proposition}

\begin{proof}
By naturality of \eqref{eq:Hopf-operation},
\eqref{eq:Hopf-on-gauge}, and the fact that $\cT$ is a group homomorphism,
\[
 \cT_*(\eta_3\eta_4)
 =\mathfrak h_{\cD}(\cT_*\eta_3).
\]
Combining \cref{lem:T-Gromoll,lem:Hopf-detection} with
\eqref{eq:exotic-eight-sphere} gives
\[
 p_9(\lambda)
 =p_9\Lambda_2\bigl(\mathfrak h_{\cD}(\cT_*\eta_3)\bigr)
 =[\epsilon]\eta=[\epsilon\eta].
\]
Stable homotopy is graded commutative, so $\epsilon\eta=\eta\epsilon$.
Relation C.2 and Table 4 of \cite{Wang} give
$[\eta\epsilon]=[\nu^3]\ne0$ in $\Coker J_9$.

The factorization \eqref{eq:alpha-cokernel-factorization} now gives
\[
 \alpha(\lambda)
 =\bar\alpha_9([\epsilon\eta])
 =\alpha_8^{\mathrm{fr}}(\epsilon)\alpha_1^{\mathrm{fr}}(\eta)=0.
\]
\end{proof}

In the Gromoll filtration, the subgroup
$\Gamma^9_3=\operatorname{im}(\Lambda_2)\subset\Theta_9$ contains
$\lambda$.  Wang proves that $p_9:\Gamma^9_3\to\Coker J_9$ is split
surjective and constructs a class with image $[\epsilon\eta]$
\cite[Corollary 3.2 and Section 4.1.2]{Wang}.

\section[The class kappa + lambda]{\texorpdfstring{The class $\kappa+\lambda$}{The class kappa + lambda}}\label{sec:last-topology}
We use \cref{lem:quotient-action-square} and quaternionic inversion to
compute $\Phi(\cA L)$.

\begin{proposition}[A relation in $\cE$]\label{prop:A2R-relation}
The class $R\in\cE$ satisfies
\begin{equation}\label{eq:A2R-local}
 \cA^2R=R\,\iota(\lambda).
\end{equation}
\end{proposition}
\begin{proof}
By \cref{lem:quotient-action-square}, $(\cA^2R)\cN=R\cN$.
Thus $\cA^2R=R\iota(w)$ for a unique $w\in\Theta_9$.
Applying $\Phi$, and using $\Phi(R)=0$ and $\Phi\iota=\id$,
identifies $w=\lambda$.
\end{proof}

\begin{lemma}[Quaternionic inversion]
\label{lem:inversion-framing}
For a smooth based map $\tau:S^3\to SO(5)$, write
$T_\tau(q,u)=(q,\tau(q)u)$.  Let $\widehat R_{\tau,b}$ and
$\widehat L_{\tau,b}$ be the diffeomorphisms of $S^8$ obtained by
transporting $R_b$ and $L_b$ through the embedding $e_0T_\tau$ and
extending by the identity outside its image.  Here
$R_b(q,u)=(qb(u),u)$ and $L_b(q,u)=(b(u)q,u)$.  Set
$\tau^\vee(q)=\tau(q^{-1})$ and
$\mathcal I(w,u)=(\overline w,u)$ on $S(\Hh\oplus\R^5)$.  Then
\begin{equation}\label{eq:inversion-exact}
 \mathcal I\widehat R_{\tau,b}\mathcal I^{-1}
       =\widehat L_{\tau^\vee,b^{-1}}.
\end{equation}
In particular, with the fixed generator $b$ and for every $m\in\Z$,
\begin{equation}\label{eq:inversion-classes}
 \Phi(\cA^mR)=\Phi(\cA^{-m}L).
\end{equation}
\end{lemma}
\begin{proof}
The ambient map $\mathcal I$ preserves the image of the standard
embedding $e_0$, since quaternionic conjugation preserves $|w|$.
In the coordinates $e_0$, it acts as $\jmath(q,u)=(q^{-1},u)$.
It therefore suffices to compute the conjugation on $T$.  Direct composition gives
\[
 \jmath T_\tau\jmath=T_{\tau^\vee},\qquad
 \jmath R_b\jmath=L_{b^{-1}}.
\]
Conjugating the product $T_\tau R_bT_\tau^{-1}$ proves
\eqref{eq:inversion-exact} on $e_0(T)$ and its complement, where both sides are the identity.
The common identity collar extends this equality smoothly to $S^8$.

The map $q\mapsto q^{-1}$ has degree $-1$ on $S^3$, hence
$[\tau^\vee]=-[\tau]$ in $\pi_3SO(5)$.  Moreover, $b^{-1}$ is homotopic
to $b$ through maps equal to $1$ near the boundary, since
$\pi_5(S^3)$ has order two.
The ambient map $\mathcal I$ reverses orientation and extends linearly
over $D^9$; conjugation by it negates the clutched sphere class.
Negation is trivial in $\Theta_9$.  Taking $[\tau]=m[\tau_1]$ gives
\eqref{eq:inversion-classes}.
\end{proof}

\begin{theorem}[Identification of $\kappa+\lambda$]\label{thm:kappa-lambda-identification}
The class $\cA L$ satisfies
\begin{equation}\label{eq:kappa-lambda-class}
 \Phi(\cA L)=\kappa+\lambda.
\end{equation}
\end{theorem}
\begin{proof}
Apply $\cA^{-1}$ to \eqref{eq:A2R-local}.  Since $\cA^{-1}$ fixes
$\cN$ pointwise, this gives
\[
 \cA R=(\cA^{-1}R)\iota(\lambda).
\]
Hence $\Phi(\cA^{-1}R)=\kappa-\lambda=\kappa+\lambda$.
Quaternionic inversion, \eqref{eq:inversion-classes} with $m=-1$,
now gives $\Phi(\cA L)=\Phi(\cA^{-1}R)$.
\end{proof}

\subsection[An attaching map for kappa + lambda]{\texorpdfstring{An attaching map for $\kappa+\lambda$}{An attaching map for kappa + lambda}}
Return to the representation $V_\kappa$ in \eqref{eq:rho-kappa}.  On the
open set $w\ne0$, set $q=w/\abs w$ and define
\begin{equation}\label{eq:b-kappa-global}
 b_\kappa(t,w,z)=b(t,q^{-1}z).
\end{equation}
If $b(t,q^{-1}z)\ne1$, then $t^2+|z|^2<c_0^2$, and hence
$|w|^2>1-c_0^2$.  Thus \eqref{eq:b-kappa-global} equals $1$ on a
neighborhood of $\{w=0\}$ and extends smoothly by that value.  It satisfies
\begin{equation}\label{eq:b-invariant}
 b_\kappa(\rho_\kappa(h)y)=b_\kappa(y).
\end{equation}

Define
\begin{equation}\label{eq:L-b}
  L_{b_\kappa}(y)=\rho_\kappa(b_\kappa(y))y.
\end{equation}
By \cref{lem:associated-diffeomorphism}, this is a diffeomorphism with
inverse $L_{b_\kappa^{-1}}$.  In the coordinates \eqref{eq:e-kappa},
\[
  e_\kappa^{-1}\circ L_{b_\kappa}\circ e_\kappa(q,t,z)
  =(b(t,z)q,t,z),
\]
so, relative to the standard embedding $e_0$, it is the class $\cA L$.
Thus \cref{thm:kappa-lambda-identification} identifies the sphere obtained by gluing
along $L_{b_\kappa}$.

\begin{corollary}\label{cor:last-clutching}
The homotopy sphere
\[
  \Sigma_{\kappa+\lambda}
  :=D^9_N\cup_{L_{b_\kappa}}D^9_S
\]
represents $\kappa+\lambda\in\Theta_9$.
\end{corollary}

\section{Positively curved metrics on the two discs}\label{sec:geometry}

We construct compatible positively curved metrics on two discs for the
conjugation-equivariant transitions in \cite{HLY,GG,DHZ}, and for the
invariant transition that gives $L_{b_\kappa}$ in
\cref{cor:last-clutching}.

For a free isometric $\star$-action on $(P,G)$, let
$\varpi:P\to P/S^3_\star$ be the quotient map.  A vector is
\emph{horizontal for $\varpi$} if it lies in $(\ker d\varpi)^{\perp_G}$.
The quotient metric $\bar G$ is defined by
\[
 \langle d\varpi_pX,d\varpi_pY\rangle_{\bar G}=\langle X,Y\rangle_G
 \quad\text{for }X,Y\in(\ker d\varpi_p)^{\perp_G}.
\]
The isometric action makes this definition independent of the point
chosen in the orbit, and $\varpi$ is a Riemannian submersion.
Its horizontal space generally differs from the horizontal space
$\ker\omega$ of the principal connection.  For a quotient plane $\Pi$ and its horizontal lift
$\widetilde\Pi$, O'Neill's formula gives
$\Sec_{\bar G}(\Pi)\geq\Sec_G(\widetilde\Pi)$ \cite{O'Neill}.

We use the following specialization of the smooth gluing theorem.
\begin{theorem}[Reiser--Wraith]\label{thm:sectional-gluing}
Let $(M_N,g_N)$ and $(M_S,g_S)$ be compact smooth Riemannian
$n$-manifolds with boundary, $n\geq3$, whose sectional curvatures are
strictly positive, including at the boundary.  Let
$\phi:\partial M_N\to\partial M_S$ be a boundary isometry, and write
$B_i(Y,Z)=\langle\nabla_Y n_i,Z\rangle_{g_i}$ for the second
fundamental form with outward unit normal $n_i$.  If
$B_N+\phi^*B_S\geq0$, then $M_N\cup_\phi M_S$ admits a smooth metric
of positive sectional curvature agreeing with $g_N,g_S$ outside any
prescribed neighborhood of the interface.
\end{theorem}
\begin{proof}
This is \cite[Theorem A(i), $k=1$]{RW}.  There $\operatorname{Ric}_k$
is a sum of $k$ sectional curvatures through a unit vector, so
$\operatorname{Ric}_1>0$ is precisely positive sectional curvature.
Their inward-normal convention
$\mathrm{II}(Y,Z)=-\langle\nabla_Y n_{\rm in},Z\rangle_{g_i}$
agrees with the outward-normal form $B_i$ above.
\end{proof}

\begin{corollary}[Equivariant gluing]\label{cor:equivariant-gluing}
Under the hypotheses of \cref{thm:sectional-gluing}, suppose a compact
Lie group acts isometrically on both pieces and $\phi$ is equivariant.
Then the glued metric can be chosen invariant.
\end{corollary}
\begin{proof}
The inward normal exponential maps give equivariant collars.  In each
collar the group acts on the boundary coordinate and fixes the distance
coordinate $t$.  Write one collar metric as $dt^2+h_t$ and choose a
smooth radial cutoff $\chi_0$ supported in that collar, equal to one
near $t=0$.  For a constant $\sigma>0$, replace its slice metrics by
\[
 h_t^{\sigma}=h_t-2\sigma t\chi_0(t)h_0.
\]
This invariant perturbation fixes the boundary metric and changes the
outward second fundamental form from $-\tfrac12h'_0$ to
$-\tfrac12h'_0+\sigma h_0$.  It tends to zero in $C^2$ as
$\sigma\to0$, with the collar and cutoff fixed.  Compactness and the
strict sectional-curvature bound therefore preserve positivity for
small $\sigma$, while the sum of boundary forms becomes positive
definite.

The cubic interpolation in \cite[Section 2.1, equation (1)]{RW} is a
linear combination, with radial scalar coefficients, of the two slice
metrics and their radial derivatives.  It is consequently invariant.
For the smoothing in \cite[Section 2.2]{RW}, use the same radial
mollifier and radial cutoff for all tensor components.  Tangential
coordinate changes are independent of $t$, so this defines smoothing
of the slice tensor itself and commutes with pullback by the group.
The curvature estimates in that construction allow common widths on
the compact boundary.  The resulting invariant smooth metric has
positive sectional curvature and agrees with the original metrics
outside the prescribed neighborhood.
\end{proof}

\subsection{Connection metrics}
We use the sign convention for $R^g$ fixed in the introduction.
In the curvature calculations, $R$ denotes the curvature tensor of the
total-space metric, and $R^B$ denotes that of the base metric.
We first state the quantitative connection-metric criterion used on the
southern disc, following \cite[Proposition 3.1]{HLY}.  The curvature
calculation is given in Appendix~\ref{app:connection-proof}.

\begin{theorem}[Positive connection metrics]\label{thm:connection-metric}
Let
\[
  \pi:P\longrightarrow(B,g_B)
\]
be a principal right $S^3$-bundle over a compact Riemannian manifold,
possibly with boundary, and assume $\Sec_{g_B}\geq\kappa_0>0$.
Put $\mathfrak g=\im\Hh$, and let $Q_{\mathfrak g}$ be its inner product for which
the induced bi-invariant metric on $S^3$ has sectional curvature one;
we use $h_{S^3}$ for that metric on the fiber.  Let $\omega$
be a smooth principal connection.  Regard its curvature $\Omega$ as an
adjoint-bundle-valued two-form on $B$; the same symbol denotes its
corresponding equivariant horizontal form on $P$.  Put
\[
  M_0=\sup_B\abs\Omega,
  \qquad
  M_1=\sup_B\abs{D\Omega}.
\]
Here $D$ uses the base Levi--Civita connection and the induced adjoint
connection.  In base and Lie-algebra orthonormal frames we use
$|\Omega|^2=\sum_{i<j,a}(\Omega^a_{ij})^2$ and
$|D\Omega|^2=\sum_{k,i<j,a}((D_k\Omega)^a_{ij})^2$.
Set $C_0=16$ and choose
\begin{equation}\label{eq:connection-concavity}
 \Lambda\geq16C_0M_0^2+
     32C_0^2\kappa_0^{-1}(M_1+1)^2+8.
\end{equation}
For $\varphi\in C^\infty(B)$, its Hessian is
$\Hess_{g_B}\varphi(X,Y)=X(Y\varphi)-(\nabla^B_XY)\varphi$.
Suppose $-\Hess_{g_B}\varphi(X,X)\geq\Lambda\langle X,X\rangle_{g_B}$ for every
base tangent vector $X$.  Then the connection metric
\begin{equation}\label{eq:connection-metric}
 G_\eps=\pi^*g_B+r_\eps^2Q_{\mathfrak g}(\omega,\omega),\qquad
 r_\eps^2=\eps e^{\eps\varphi}
\end{equation}
has strictly positive sectional curvature whenever
$0<\eps\leq\eps_{\mathrm{conn}}$, where the positive number
$\eps_{\mathrm{conn}}$ is given explicitly in
\eqref{eq:conn-epsilon-threshold}.  More precisely,
\begin{equation}\label{eq:connection-margin-main}
 \Sec_{G_\eps}\geq
 \min\left\{\frac{\kappa_0}{2},\frac{\eps\Lambda}{8},
                         \frac1{4\eps}\right\}>0.
\end{equation}
The functions in \eqref{eq:connection-metric} are pulled back from $B$.
\end{theorem}

\begin{proof}
The curvature computation and simultaneous estimates in
Appendix~\ref{app:connection-proof} give
\eqref{eq:conn-final-area-lower} under the stated hypotheses.
To deduce the sectional bound, take independent tangent vectors
$P=X+U$ and $Q=Y+V$, split into horizontal and vertical parts with respect to $\omega$.
With norms induced by $G_\eps$, put
\[
 h=|X\wedge Y|,\qquad
 m=|X\otimes V-Y\otimes U|,\qquad k=|U\wedge V|.
\]
These are the norms of the three orthogonal components of $P\wedge Q$,
so $|P\wedge Q|^2=h^2+m^2+k^2$.  The cited estimate states
\[
 \langle R(P,Q)Q,P\rangle_{G_\eps}
 \geq\frac{\kappa_0}{2}h^2+\frac{\eps\Lambda}{8}m^2
          +\frac{1}{2r_\eps^2}k^2.
\]
The same scale choice gives $r_\eps^2\leq2\eps$.  Dividing by the
positive squared area $|P\wedge Q|^2$ proves
\eqref{eq:connection-margin-main}.
\end{proof}

\subsection[Two actions on the bundle]{\texorpdfstring{Two $S^3$-actions on the bundle}{Two actions on the bundle}}\label{subsec:bundle-actions}
Let $V$ be an orthogonal real representation of $S^3$ of dimension
$n\geq3$, and let
\[
  \rho:S^3\longrightarrow O(V).
\]
On the unit sphere $S(V)$ define the infinitesimal action map
\begin{equation}\label{eq:K-def}
  K_y:\im\Hh\longrightarrow T_yS(V),
  \qquad
  K_y\xi=\left.\frac{\dd}{\dd s}\right|_{s=0}
  \rho(e^{s\xi})y.
\end{equation}
We use the unit round metric on $S(V)$ and the standard Euclidean metric
on $\im\Hh$, and fix
\begin{equation}\label{eq:K-bound}
 L_\rho=\sup_{y\in S(V)}\norm{K_y}_{\op}<\infty.
\end{equation}

Let $U_N,U_S$ be the closed northern and southern hemispheres of
$S(\R\oplus V)$, with common equator $S(V)$.  The base action is trivial
on $\R$ and is $\rho$ on $V$.  We use the equivariant disc charts
\[
 x\longmapsto\left(\pm\cos\frac{\pi|x|}{2},\,
              \frac{\sin(\pi|x|/2)}{|x|}x\right),\qquad x\in D(V),
\]
with the second factor interpreted smoothly at $x=0$.
They agree on the equator; $x$ below denotes these disc coordinates.
For a smooth map
$\gamma:S(V)\to S^3$, define the principal right bundle $P_\gamma$ by
gluing $U_N\times S^3$ to $U_S\times S^3$ with the identification
\[
 (y,v_N)_N\sim(y,\gamma(y)v_N)_S.
\]
The transition is extended independently of the signed radial variable
on a collar of the equator; this specifies the smooth bundle structure.
The right action is $(x,v_i)h=(x,v_ih)$ in either chart.

For a conjugation-equivariant transition $\gamma=\theta$, the local
formula
\[
 q\star(x,v_i)=(\rho(q)x,qv_i)
\]
defines a global free left action on $P_\theta$.  Indeed,
$\theta(\rho(q)y)qv_N=q\theta(y)v_N$, which is the compatibility
condition on the overlap.  It commutes with the principal right action.
The quotient coordinate is $\rho(v_i^{-1})x$.  On the equator the change
of quotient coordinates is therefore
$J_{\theta^{-1}}(y)=\rho(\theta(y)^{-1})y$.

For an invariant transition $\gamma=b$, use instead
\begin{equation}\label{eq:right-diagonal-star}
 q\star(x,v_i)=(\rho(q)x,v_iq^{-1}).
\end{equation}
Invariance of $b$ makes this formula compatible with the transition.
It is a free left action, since $v_iq^{-1}=v_i$ forces $q=1$.
Its quotient coordinate is $\rho(v_i)x$, and the change of quotient
coordinates on the equator is $L_b(y)=\rho(b(y))y$.
This action need not commute with the principal right action.
The same formulas hold in polar coordinates $x=sy$ and extend to $x=0$.

\subsection{Quotient metrics and boundary coordinates}
We compute the quotient metric on a fixed radial level.

\begin{lemma}[The quotient metric]\label{lem:quotient-metric}
Let $f,r>0$.  On $S(V)\times S^3$ use the product metric
$f^2h_{S(V)}+r^2h_{S^3}$ and either diagonal action described above.
Use the section $y\mapsto(y,1)$ to identify the quotient with $S(V)$.
Write
\[
 D_y=r^2I+f^2K_y^*K_y,\qquad
 C_y=I+\frac{f^2}{r^2}K_yK_y^*.
\]
For $Y,Z\in T_yS(V)$ the quotient metric is
\begin{align}\label{eq:quotient-metric}
 h_{f,r}(Y,Z)
 &=f^2\ip{C_y^{-1}Y}{Z}\notag\\
 &=f^2\ip{Y}{Z}
   -f^4\ip{D_y^{-1}K_y^*Y}{K_y^*Z}.
\end{align}
In particular,
\begin{equation}\label{eq:boundary-round-comparison}
 \frac{f^2r^2}{r^2+f^2L_\rho^2}\,h_{S(V)}
 \ \leq\ h_{f,r}\ \leq\ f^2h_{S(V)}.
\end{equation}
The angular component of the horizontal lift of $Y$, measured in the
scaled orthonormal frame, is $X=fC_y^{-1}Y$.
\end{lemma}
\begin{proof}
Fiber inversion identifies the two actions isometrically and fixes the
section, so it suffices to use the diagonal action $q\star(x,v)=(\rho(q)x,qv)$.
Project the section vector $(Y,0)$ orthogonally away from the orbit
$(K_y\xi,\xi)$.  Orthogonality gives
\[
 \xi=f^2D_y^{-1}K_y^*Y,
 \qquad (Y,0)^{\mathrm{hor}}=(Y-K_y\xi,-\xi).
\]
The identity
\[
 I-f^2K_yD_y^{-1}K_y^*=C_y^{-1}
\]
follows by multiplication by $C_y$.  Taking the squared norm of this
horizontal lift, and then polarizing, proves
\eqref{eq:quotient-metric}.
Since $I\leq C_y\leq(1+f^2L_\rho^2/r^2)I$, it also proves
\eqref{eq:boundary-round-comparison}.  Since $D_y\geq r^2I$ and
$C_y\geq I$, their inverses depend smoothly on $y$.
\end{proof}

\begin{lemma}[The metric on a boundary collar]\label{lem:collar-metric}
Let $\gamma$ be either type of transition in the preceding construction,
and let $t$ be a collar coordinate.  Suppose $\gamma$ is independent of
$t$ and the northern-gauge metric is
$dt^2+f(t)^2h_{S(V)}+r(t)^2h_{S^3}$, where $f,r>0$.  In the southern gauge
$v_S=\gamma(y)v_N$, the connection potential is
$A_S=-d\gamma\,\gamma^{-1}$.  Put $\phi=J_{\gamma^{-1}}$ in the
conjugation-equivariant case and $\phi=L_\gamma$ in the invariant case.
Then the quotient metric in the coordinates extending over the southern disc is
\begin{equation}\label{eq:collar-pullback}
 \bar G_S=dt^2+(\phi^{-1})^*h_{f(t),r(t)}.
\end{equation}
\end{lemma}
\begin{proof}
The change of gauge
$F(t,y,v_N)=(t,y,\gamma(y)v_N)$ satisfies
\[
 F^*\bigl(\Ad_{v_S^{-1}}A_S+v_S^{-1}dv_S\bigr)
 =v_N^{-1}dv_N.
\]
Thus the change of trivialization is an isometry between the two
expressions for the total-space metric.  It is equivariant for the corresponding star action and,
in quotient coordinates, becomes $(t,y)\mapsto(t,\phi(y))$.
The quotient metric in \cref{lem:quotient-metric} therefore gives
\eqref{eq:collar-pullback}.  Since $\phi$ is independent of $t$,
radial differentiation commutes with this pullback.
\end{proof}

\subsection{Curvature of the northern quotient}
The northern metric is a doubly warped product.  The following criterion
controls its curvature on planes horizontal for $\varpi$.

\begin{proposition}[Positive curvature of the northern quotient]\label{prop:north-criterion}
Let $\ell>0$ and consider
\begin{equation}\label{eq:north-metric}
 G_N=ds^2+F(s)^2h_{S(V)}+r(s)^2h_{S^3},\qquad 0\leq s\leq\ell,
\end{equation}
where $F,r\in C^\infty([0,\ell])$, $F(0)=0$, $F'(0)=1$, and $F$ has
an odd smooth extension near zero.  Assume $r>0$ is constant near zero,
$F'>0$, and $0\leq r'<1$.  For a real number $a$, write
$a_+=\max\{a,0\}$.
On $s>0$ let $S^3$ act freely by
\begin{equation}\label{eq:inverted-star}
 q\star(y,u)=(\rho(q)y,qu).
\end{equation}
Suppose the two functions
\begin{equation}\label{eq:north-sufficient}
 \begin{split}
 c_{\mathrm{rad}}&=-\frac{F''}{F}-\frac{L_\rho^2F^2(r'')_+}{r^3}>0,\\
 c_{\mathrm{ang}}&=\frac{1-F'^2}{F^2}-\frac{2L_\rho^2FF'r'}{r^3}>0
 \end{split}
\end{equation}
extend continuously and positively to $s=0$.  Then $G_N$ extends
smoothly across the center and the quotient is a smooth positively
curved disc.  Set $F_a=F(\ell)$,
$r_{\min}=\min r$, and
$m_N=\min_{[0,\ell]}\min\{c_{\mathrm{rad}},c_{\mathrm{ang}}\}$.
Its sectional curvature is at least
\[
 \frac{m_N}{(1+L_\rho^2F_a^2/r_{\min}^2)^2}>0.
\]
\end{proposition}
\begin{proof}
In Cartesian polar coordinates $x=sy$, the quotient is identified with
$D(V)$ by $[(x,u)]\mapsto\rho(u^{-1})x$, after a fixed radial rescaling.
The global section is $x\mapsto(x,1)$.  Right multiplication
on the fiber commutes with \eqref{eq:inverted-star}, so it suffices to
compute at fiber coordinate $u=1$.  Write angular and fiber components in
orthonormal frames for the \emph{scaled} metrics $F^2h_{S(V)}$ and
$r^2h_{S^3}$.  The orbit vector associated with $\xi$ has components
$(FK_y\xi,r\xi)$ in these frames.  A vector horizontal for $\varpi$ consequently
has the form
\begin{equation}\label{eq:horizontal-graph}
  x_0\partial_s+X+T_yX,
  \qquad T_y=-\frac{F}{r}K_y^*,
  \qquad \norm{T_y}_{\op}\leq\frac{L_\rho F}{r},
\end{equation}
where $L_\rho$ is fixed as in \eqref{eq:K-bound}.

Writing the infinitesimal action as $\mathcal O_y\xi=(FK_y\xi,r\xi)$,
the orthogonal projection onto the horizontal space is
\begin{equation}\label{eq:orbit-projector}
 \mathcal O_y^*\mathcal O_y=r^2I+F^2K_y^*K_y\geq r^2I,
 \qquad
 \Pi_y=I-\mathcal O_y(\mathcal O_y^*\mathcal O_y)^{-1}\mathcal O_y^*.
\end{equation}
The lower bound makes $\Pi_y$ smooth.  The radial direction is horizontal.

Fix $s>0$ and $y\in S(V)$, and write $T=T_y$.
Take horizontal vectors
$P=x_0\partial_s+X+TX$ and $Q=y_0\partial_s+Y+TY$.
Put $A=x_0Y-y_0X$ and $B=X\wedge Y$.  The radial--angular,
angular--angular, radial--fiber, angular--fiber and fiber--fiber
components of $P\wedge Q$ are, respectively,
$A$, $B$, $TA$, $X\otimes TY-Y\otimes TX$, and $TX\wedge TY$.
The curvature operator is diagonal on these five blocks, so
\begin{equation}\label{eq:north-five-blocks}
\begin{split}
 \langle R^{G_N}(P,Q)Q,P\rangle_{G_N}={}&
 -\frac{F''}{F}|A|^2+\frac{1-F'^2}{F^2}|B|^2
 -\frac{r''}{r}|TA|^2\\
 &-\frac{F'r'}{Fr}|X\otimes TY-Y\otimes TX|^2
 +\frac{1-r'^2}{r^2}|TX\wedge TY|^2.
\end{split}
\end{equation}
The mixed component satisfies
\[
 |X\otimes TY-Y\otimes TX|^2\leq2\norm{T}_{\op}^2|X\wedge Y|^2.
\]
If $X,Y$ are dependent, both sides vanish.  Otherwise replace them by
an orthonormal basis $e_1,e_2$ of their span and divide by
$|X\wedge Y|^2$.  The left side becomes
$|Te_1|^2+|Te_2|^2\leq2\|T\|_{\op}^2$.
Write $\mathcal K(P,Q)=\langle R^{G_N}(P,Q)Q,P\rangle_{G_N}$.
The last term of \eqref{eq:north-five-blocks} is nonnegative because
$0\leq r'<1$.  The radial--fiber term is bounded below by
$-(r'')_+\|T\|_{\op}^2|A|^2/r$.  Since $F'>0$ and $r'\geq0$, the
angular--fiber coefficient is nonpositive, so the preceding upper bound
on its area gives a lower bound for its curvature contribution.  Using
$\|T\|_{\op}\leq L_\rho F/r$ in these two estimates yields
\begingroup
\allowdisplaybreaks[0]
\begin{align}\label{eq:north-numerator}
 \mathcal K(P,Q)\geq{}&
 \left(-\frac{F''}{F}-\frac{L_\rho^2F^2(r'')_+}{r^3}\right)|A|^2\notag\\
 &+\left(\frac{1-F'^2}{F^2}
             -\frac{2L_\rho^2FF'r'}{r^3}\right)|B|^2.
\end{align}
\endgroup
The projected bivector is
\[
 (x_0\partial_s+X)\wedge(y_0\partial_s+Y)
 =\partial_s\wedge A+B.
\]
If $A=B=0$, the two projected vectors are dependent.  Since the
horizontal space is the graph of a linear map over this projection,
$P,Q$ would also be dependent.  Thus every horizontal two-plane has a
positive numerator under the stated hypotheses.  At $s=0$, the metric $G_N$ is smooth because $F$ has an odd smooth extension, $F'(0)=1$,
and $r$ is constant nearby.  In Cartesian disc coordinates, the orbit at the center is contained
in the fiber, so the horizontal space is the base tangent space.
The sectional curvature of $G_N$ on such a plane is the common positive limit $-F'''(0)$ of the two base curvature
coefficients.

With $m_N$ and $r_{\min}$ as in the statement, the two coefficients
in \eqref{eq:north-numerator} are at least $m_N>0$.
The map $(x_0,X)\mapsto(x_0,X,TX)$ has norm at most
$\sqrt{1+\norm T_{\op}^2}$.  Taking its exterior square gives
\begin{equation}\label{eq:graph-area}
 |P\wedge Q|^2\leq
 (1+L_\rho^2F_a^2/r_{\min}^2)^2\bigl(|A|^2+|B|^2\bigr).
\end{equation}
Consequently, on every plane horizontal for $\varpi$, including at the
center, the sectional curvature of $G_N$ is at least
\begin{equation}\label{eq:north-uniform-margin}
 c_N=\frac{m_N}{(1+L_\rho^2F_a^2/r_{\min}^2)^2}>0.
\end{equation}
O'Neill's formula \cite{O'Neill} gives the same lower bound on the northern
quotient.

\end{proof}

Appendix~\ref{app:warping-functions} constructs $F$ and $r$ with prescribed
boundary radii and logarithmic derivative $r^\prime/r$.

\subsection{Construction and gluing of the disc metrics}
For either transition function, choose the latitude $a$, cutoff $\chi$
and principal connection $\omega$ first.  Their curvature bounds fix
$M_0,M_1$ and then $\Lambda,A_0$.  Next choose the cutoff $\zeta$ and
$\delta,F,\ell$ in Appendix~\ref{app:warping-functions}; the common
fiber scale $\eps$ is chosen last, using both the connection bound and
$\ell$.  All these choices may depend on the fixed transition and
representation.

\begin{theorem}[Invariant transition functions]\label{thm:invariant-transition}
Let $b:S(V)\to S^3$ be smooth and $\rho$-invariant.  Let
$P_b\to S(\R\oplus V)$ be the principal right $S^3$-bundle with
equatorial transition $v_S=b(y)v_N$, equipped with the action \eqref{eq:right-diagonal-star}.  Its quotient
admits a smooth Riemannian metric with strictly positive sectional
curvature.  With the disc parametrizations specified above, it is
\begin{equation}\label{eq:invariant-quotient}
  P_b/S^3_\star\cong D^n_N\cup_{L_b}D^n_S,
  \qquad n=\dim V.
\end{equation}
\end{theorem}

\begin{proof}
\smallskip
\noindent\emph{The quotient and the attaching map.}
For the action of Section~\ref{subsec:bundle-actions}, the local quotient of $D(V)\times S^3$ is the smooth disc $D(V)$,
with quotient map $(x,v)\mapsto\rho(v)x$ and section $x\mapsto(x,1)$.
In either polar trivialization define the quotient
angular coordinate
\begin{equation}\label{eq:quotient-coordinate}
  y_i=\rho(v_i)y.
\end{equation}
It is invariant under $\star$.  On the overlap,
\[
  y_S=\rho(b(y)v_N)y
      =\rho(b(y_N))y_N
      =L_b(y_N),
\]
where the middle equality uses the $\rho$-invariance of $b$.  This proves
\eqref{eq:invariant-quotient}.

\smallskip
\noindent\emph{The southern connection.}
Fix
\[
  \frac\pi2<a<b_0<\pi.
\]
Transport the equatorial transition radially to the latitude $t=a$.
The two base pieces $0\leq t\leq a$ and $a\leq t\leq\pi$ are still
standard discs with the transported parametrizations.  Radial reparameterization does not change their
angular clutching map.  Use the round southern base disc
\begin{equation}\label{eq:south-base}
  g_B=\dd t^2+\sin^2t\,h_{S(V)},
  \qquad a\leq t\leq\pi.
\end{equation}
Choose a smooth cutoff $\chi:[a,\pi]\to[0,1]$ that is zero near $t=a$
and one for $t\geq b_0$.
On the punctured cap $a\leq t<\pi$, pull $b$ back independently of $t$
and put, in the northern and southern gauges,
\begin{equation}\label{eq:connection-potentials}
  A_N=\chi b^{-1}\dd b,
  \qquad
  A_S=(\chi-1)\dd b\,b^{-1}.
\end{equation}
Since $\chi=1$ near $\pi$, extend $A_S$ by zero across the south pole.
The relation
\[
  A_N=b^{-1}A_Sb+b^{-1}\dd b
\]
shows that these potentials define a smooth principal connection: in
the $i$th trivialization its connection form is
$\omega_i=\Ad_{v_i^{-1}}A_i+v_i^{-1}dv_i$.
It is a product connection near the boundary in the northern gauge
and near the center in the southern gauge.  Put $a_b=b^{-1}db$,
pulled back independently of $t$, and use the convention
$(a_b\wedge a_b)(Y,Z)=[a_b(Y),a_b(Z)]$.  The Maurer--Cartan identity
$da_b+a_b\wedge a_b=0$ gives the local curvature forms
\begin{equation}\label{eq:cutoff-connection-curvature}
 \Omega_N=\chi'\,dt\wedge a_b+\chi(\chi-1)a_b\wedge a_b,
 \qquad \Omega_S=\operatorname{Ad}_b\Omega_N.
\end{equation}
The curvature vanishes where $\chi$ is identically $0$ or $1$.
Thus $\Omega$ and $D\Omega$ are supported in the fixed annulus
$[a,b_0]\times S(V)$, where $\sin t$ is bounded away from zero.
Their norms, taken with respect to the fixed round base and Lie-algebra
metrics, are bounded independently of $r$ and $\eps$.

Let $B_q(y,v)=(\rho(q)y,v)$ and let $R_h(y,v)=(y,vh)$ denote the
principal right action.  The invariance of $b$ makes $B_q$ compatible
with the transition and gives $B_q^*A_i=A_i$ in both gauges.
Thus $B_q^*\omega=\omega$, whereas
$R_h^*\omega=\Ad_{h^{-1}}\omega$.  Since the $\star$-action is
$B_qR_{q^{-1}}$, it follows that
\begin{equation}\label{eq:omega-star}
  (q\star)^*\omega=\Ad_q\omega,
  \qquad
  (q\star)^*Q_{\mathfrak g}(\omega,\omega)=Q_{\mathfrak g}(\omega,\omega).
\end{equation}
Let $M_0,M_1$ be the finite curvature bounds of this fixed connection.
Choose $\Lambda$ as in \cref{thm:connection-metric} and choose $A_0>0$
with
\[
  A_0\abs{\cos a}\geq\Lambda.
\]
The round height function satisfies $\Hess(\cos t)=-\cos t\,g_B$.
Thus $\varphi(t)=-A_0\cos t$ satisfies
\[
 -\Hess\varphi=-A_0\cos t\,g_B
 \geq A_0|\cos a|\,g_B\geq\Lambda g_B,
 \qquad a\leq t\leq\pi.
\]
Theorem~\ref{thm:connection-metric} therefore makes
\begin{equation}\label{eq:south-metric}
 G_S=\dd t^2+\sin^2t\,h_{S(V)}+r_S(t)^2Q_{\mathfrak g}(\omega,\omega),
 \qquad
 r_S(t)^2=\eps e^{-\eps A_0\cos t},
\end{equation}
a metric of strictly positive sectional curvature for every
$0<\eps\leq\eps_{\mathrm{conn}}$.  The action is free by its definition;
\eqref{eq:omega-star} and invariance of the round base and $r_S$ show
that it is isometric.  O'Neill's formula gives positive sectional
curvature on the southern quotient disc.

\smallskip
\noindent\emph{The northern metric.}
Write $r_S(t)=\sqrt\eps\,e^{-\eps A_0\cos t/2}$ and choose
$F,r,\ell$ by Appendix~\ref{app:warping-functions} with boundary data
$F_a=\sin a$, $r_a=r_S(a)$, and $q_s=r'_S(a)/r_a$.  Fiber inversion $u=v^{-1}$ is an isometry of the
round fiber and converts the action \eqref{eq:right-diagonal-star} to
\eqref{eq:inverted-star} while fixing the section $x\mapsto(x,1)$.
\Cref{prop:north-criterion} therefore gives a smooth strictly positively
curved quotient disc, including its center.  These choices also give
$F(\ell)=F_a$, $r(\ell)=r_a$, and $r'(\ell)/r_a=q_s$.

\smallskip
\noindent\emph{Boundary data and gluing.}
In the northern trivialization, both total-space boundary metrics are
\[
  F_a^2h_{S(V)}+r_a^2h_{S^3},
  \qquad F_a=\sin a,
\]
and the $\star$-action is the same.  In the coordinates extending over the discs, the
quotient boundary metrics therefore satisfy
\begin{equation}\label{eq:boundary-isometry}
  h_N=L_b^*h_S.
\end{equation}
On a southern collar, \cref{lem:collar-metric} gives
\[
 \bar G_S=dt^2+(L_b^{-1})^*h_{\sin t,r_S(t)}.
\]
On the northern collar the quotient is $ds^2+h_{F(s),r(s)}$.
The angular identification $L_b$ is independent of both radial coordinates.
Put
\[
 \mu_N=\frac{F'(\ell)}{F_a}>0,\qquad
 \mu_S=\cot a<0,\qquad
 q_s=\frac{r'(\ell)}{r_a}=\frac{r_S'(a)}{r_a}.
\]
In the common boundary coordinates, write
$h=F^2C^{-1}$, where $C=I+(F^2/r^2)K_yK_y^*$.
For $\mu=F'/F$ and $q=r'/r$,
$C'=2(\mu-q)(C-I)$ commutes with $C$.  Hence
\begin{equation}\label{eq:individual-shape}
 \tfrac12h^{-1}h'=qI+(\mu-q)C^{-1}.
\end{equation}
The outward normals are $+\partial_s$ on the northern disc and
$-\partial_t$ on the southern disc.  Their shape endomorphisms are
$q_sI+(\mu_N-q_s)C_y^{-1}$ and
$-q_sI-(\mu_S-q_s)C_y^{-1}$, respectively, where
$C_y=I+(F_a^2/r_a^2)K_yK_y^*$.  Their sum is
\begin{equation}\label{eq:boundary-shape-endomorphism}
 h_N^{-1}\bigl(B_N+L_b^*B_S\bigr)
 =(\mu_N-\mu_S)C_y^{-1}.
\end{equation}
Since $C_y\leq(1+L_\rho^2F_a^2/r_a^2)I$, we obtain
\begin{equation}\label{eq:boundary-uniform-margin}
 B_N+L_b^*B_S\ \geq\
 \frac{\mu_N-\mu_S}{1+L_\rho^2F_a^2/r_a^2}\,h_N>0.
\end{equation}
Together with \eqref{eq:boundary-isometry}, this verifies the hypotheses
of \cref{thm:sectional-gluing} and proves the result.
\end{proof}

\begin{theorem}[Conjugation-equivariant transition functions]\label{thm:equivariant-transition}
Let $\rho:S^3\to O(V)$ be an orthogonal representation with
$\dim V=n\geq3$, and let $\theta:S(V)\to S^3$ be smooth and
conjugation equivariant.  Form $P_\theta\to S(\R\oplus V)$ with
transition $v_S=\theta(y)v_N$ and the action $q\star(x,v)=(\rho(q)x,qv)$.
Then the quotient
\[
 P_\theta/S^3_\star\cong D_N^n\cup_{J_{\theta^{-1}}}D_S^n
\]
admits a smooth metric of strictly positive sectional curvature.
\end{theorem}

\begin{proof}
This construction also appears in \cite[Proposition 5.1]{DHZ} and
\cite[Theorem 4.1 and Remark 4.4]{GG}.
Conjugation equivariance makes $v_S=\theta(y)v_N$ compatible with
the action $q\star(x,v)=(\rho(q)x,qv)$.  Its quotient disc
coordinate is $y_i=\rho(v_i^{-1})y$, including the analogous Cartesian
coordinate at the center.  Conjugation equivariance gives
\[
 y_S=\rho(\theta(y_N)^{-1})y_N=J_{\theta^{-1}}(y_N).
\]
Choose $\pi/2<a<b_0<\pi$ and a smooth cutoff
$\chi:[a,\pi]\to[0,1]$ that is zero near $a$ and one for $t\geq b_0$.
On the punctured cap take $A_N=\chi\theta^{-1}d\theta$ and
$A_S=(\chi-1)d\theta\,\theta^{-1}$, extending $A_S$ by zero over the
pole.  Under the base action their pullbacks are $\Ad_q A_i$;
simultaneous left multiplication of $v_i$ by $q$ then
cancels this adjoint factor in
$\omega_i=\Ad_{v_i^{-1}}A_i+v_i^{-1}dv_i$.  Hence the principal connection,
and therefore its connection metric, is invariant under $\star$.  Formula
\eqref{eq:cutoff-connection-curvature} applies with $b$ replaced by
$\theta$, so its curvature bounds are again fixed before the fiber scale.

Let $M_0=\sup|\Omega|$ and $M_1=\sup|D\Omega|$ for this connection.
Choose $\Lambda$ by \eqref{eq:connection-concavity}, and then choose $A_0$ with
$A_0|\cos a|\geq\Lambda$.  The southern connection metric is positive by
\cref{thm:connection-metric}.  With these data fixed, choose
$\delta,F,\ell$ and finally the common scale $\eps$ as in
Appendix~\ref{app:warping-functions}.
The northern action is already \eqref{eq:inverted-star}, so
\cref{prop:north-criterion} gives the northern metric.  The endpoint
conditions make both total-space boundary metrics equal in the northern
trivialization.  Thus $h_N=J_{\theta^{-1}}^*h_S$ on the boundary.
The southern collar metric is
$dt^2+(J_{\theta^{-1}}^{-1})^*h_{\sin t,r_S(t)}$ by
\cref{lem:collar-metric}.  Applying
\eqref{eq:individual-shape} with the two outward normals gives
\[
 h_N^{-1}(B_N+J_{\theta^{-1}}^*B_S)
   =(\mu_N-\mu_S)C_y^{-1}>0.
\]
Theorem~\ref{thm:sectional-gluing} now gives a smooth positively curved
metric on the stated clutched manifold.
\end{proof}

\begin{remark}[The induced group action]\label{rem:residual-action}
For a conjugation-equivariant $\theta$, the principal right action
commutes with $\star$ and induces a left action on the quotient by
$h\cdot[p]=[ph^{-1}]$.  In each disc coordinate this is $y\mapsto\rho(h)y$;
equivalently, $J_{\theta^{-1}}$ is $\rho$-equivariant.
The source metrics are principal-right-invariant, so both quotient
disc metrics are invariant under this action.  Corollary~\ref{cor:equivariant-gluing}
gives an invariant metric of positive sectional curvature.
If $\rho$ is nontrivial, its effective image is $S^3$ or $SO(3)$.
The Lawson--Yau theorem applied to this nonabelian group action gives
a metric of positive scalar curvature \cite{LawsonYau}.  Hence a spin
quotient has vanishing Hitchin invariant.  If $\rho$ is trivial, the
attaching map is the identity and the quotient is the standard sphere.
\end{remark}

\section{Completion of the proof}\label{sec:applications}\label{sec:proof-main}

We apply the preceding metric constructions to the three attaching maps.  The required
representation bounds are as follows.

\subsection{The two representations}
For $V_\kappa=\R\oplus\Hh^2$ and $y=(t,w,z)$, differentiation gives
\begin{equation}\label{eq:K-kappa}
 K_y\xi=(0,\xi w,\xi z).
\end{equation}
Multiplicativity of the quaternionic norm yields
\begin{equation}\label{eq:Gram-kappa}
 |K_y\xi|^2=(|w|^2+|z|^2)|\xi|^2=(1-t^2)|\xi|^2,
\end{equation}
and hence
\begin{equation}\label{eq:bound-kappa}
 K_y^*K_y=(1-t^2)I,\qquad\norm{K_y}_{\op}\leq1.
\end{equation}
For $V_\lambda=\R^2\oplus\Hh\oplus\im\Hh$ and $y=(a,w,z)$,
\begin{equation}\label{eq:K-lambda}
 K_y\xi=(0,\xi w,[\xi,z]).
\end{equation}
Since $[\xi,z]=2\xi\times z$, one obtains
\begin{equation}\label{eq:Gram-lambda}
 K_y^*K_y=|w|^2I+4\bigl(|z|^2I-z\otimes z\bigr).
\end{equation}
If $z\ne0$, the decomposition
$\im\Hh=\R z\oplus z^\perp$ gives the eigenvalue $|w|^2$ on
$\R z$ and $|w|^2+4|z|^2$ on the two-dimensional space $z^\perp$.
If $z=0$, the Gram matrix is $|w|^2I$, so its only eigenvalue has
multiplicity three.  In both cases,
$\norm{K_y}_{\op}^2=|w|^2+4|z|^2\leq4$, since
$|a|^2+|w|^2+|z|^2=1$.  Thus
\begin{equation}\label{eq:bound-lambda}
 \norm{K_y}_{\op}\leq2.
\end{equation}
The ranks are three for $w\ne0$, two for $w=0$, $z\ne0$, and zero
for $w=z=0$.

\subsection{Application to the three sphere classes}
Table~\ref{tab:attaching-maps} lists the transition functions and the
attaching maps in the specified disc coordinates.  Its column headed
``transition'' records $\gamma$ in $v_S=\gamma(y)v_N$.
The first two transitions are conjugation equivariant and use the
action $q\star(x,v)=(\rho(q)x,qv)$.  The last is invariant and uses the
action \eqref{eq:right-diagonal-star}.  The quotient-coordinate convention accounts for the inverses in the
first two transitions.

\begin{table}[ht]
\centering\small
\renewcommand{\arraystretch}{1.18}
\begin{tabular}{@{}c c c c@{}}
\toprule
sphere class & representation & transition & attaching map\\
\midrule
$\kappa$ & $V_\kappa$ & $\beta_\kappa^{-1}$ & $J_{\beta_\kappa}$\\
$\lambda$ & $V_\lambda$ & $\beta_\lambda^{-1}$ & $J_{\beta_\lambda}$\\
$\kappa+\lambda$ & $V_\kappa$ & $b_\kappa$ & $L_{b_\kappa}$\\
\bottomrule
\end{tabular}
\caption{Transition functions and attaching maps for the three nonstandard classes.}
\label{tab:attaching-maps}
\end{table}

\begin{proof}[Proof of \cref{thm:main-intro}]
The standard sphere has its round metric.  Apply \cref{thm:equivariant-transition}
to the first two rows of Table~\ref{tab:attaching-maps}; their smooth
types are identified in \cref{prop:kervaire-identification,prop:lambda-invariants}.
Apply \cref{thm:invariant-transition} to the last row, whose type is given by
\cref{cor:last-clutching}.  Thus every class of $\ker\alpha$ admits
positive sectional curvature.

Positive sectional curvature implies nonnegative sectional curvature
and positive scalar curvature.  Conversely, Hitchin's theorem gives
$\alpha=0$ from positive scalar curvature, while
\cref{prop:no-nonnegative} gives $\alpha=0$ from nonnegative sectional
curvature.  These implications prove all four equivalences.
\end{proof}

\begin{remark}\label{rem:symmetry}
On the Kervaire sphere, the metric constructed here cannot be invariant
under the full natural $SO(2)\times SO(5)$ action on its Brieskorn model.
Grove--Verdiani--Wilking--Ziller exclude nonnegative sectional curvature
for metrics invariant under that action, not for arbitrary metrics on
the underlying manifold \cite{GVWZ}.
\end{remark}

\appendix

\section{A quantitative connection-metric criterion}\label{app:connection-proof}

We prove \cref{thm:connection-metric} using the curvature formulas and
decomposable-bivector estimates of \cite[Section 3]{HLY}.  Throughout this
appendix, write $G=G_\eps$ and $r=r_\eps$.  Horizontal vectors are
identified with their projections when inserted into base tensors.
The dual metric $g_B^{-1}$ pairs base covectors.

\subsection{Frames and curvature blocks}
Write the base dimension as $d$ and let $X_i$ be the horizontal lifts, with respect to $\omega$, of a base
orthonormal frame normal at the point under consideration.
Use the fundamental vector fields $\xi^\#$ of the principal right
action defined in the bundle conventions.
Choose a positively oriented orthonormal basis $e_1,e_2,e_3$ of
$\im\Hh$ and set
$E_a=r^{-1}e_a^\#$, so that $E_1,E_2,E_3$ are unit vertical vectors.
Write
\[
 [e_a,e_b]=c_{ab}^{\ \ c}e_c,\qquad
 c_{ab}^{\ \ c}=2\epsilon_{abc},\qquad
 \vartheta=d\log r,\qquad \mathsfH=-r^{-1}\Hess r.
\]
Here $\epsilon_{abc}$ is the alternating symbol, with $\epsilon_{123}=1$;
$\vartheta^\sharp$ denotes the $g_B$-dual vector field of $\vartheta$.
Superscripts $\mathrm H,\mathrm V$ denote the orthogonal projections to
$\ker\omega$ and $\ker d\pi$.  Write $\Omega_{ij}=\sum_b\Omega_{ij}^{\,b}e_b$ in the fixed Lie-algebra
basis.  Principal equivariance gives the vertical derivatives of these
coefficient functions:
\begin{equation}\label{eq:conn-vertical-derivative}
 E_a(\Omega_{ij}^{\,b})
 =-r^{-1}\sum_c c_{ac}^{\ \ b}\Omega_{ij}^{\,c}.
\end{equation}
The bracket and Koszul identities are
\begin{align}\label{eq:connection-koszul}
 [X_i,X_j]^{\mathrm V}&=-r\Omega_{ij}^{\,a}E_a,
 & [X_i,E_a]&=-\vartheta_iE_a,\notag\\
 [E_a,E_b]&=r^{-1}c_{ab}^{\ \ c}E_c,
 & \nabla_{X_i}X_j&=(\nabla^B_iX_j)^{\mathrm H}
                    -\tfrac r2\Omega_{ij}^{\,a}E_a,\notag\\
 \nabla_{X_i}E_a&=\tfrac r2\Omega_{ij}^{\,a}X_j,
 &\nabla_{E_a}X_i&=\tfrac r2\Omega_{ij}^{\,a}X_j+\vartheta_iE_a,\\
 \nabla_{E_a}E_b&=\tfrac1{2r}c_{ab}^{\ \ c}E_c
                    -\delta_{ab}\vartheta_iX_i.\notag
\end{align}
Repeated indices are summed over their full ranges.  Put
$R_{ABCD}=\langle R(A,B)C,D\rangle_G$.  Differentiating
\eqref{eq:connection-koszul}, using
\eqref{eq:conn-vertical-derivative}, gives
\begin{align}\label{eq:connection-curvature-components}
 R_{i a b j}
 &=\mathsfH_{ij}\delta_{ab}
   +\frac{r^2}{4}\sum_l\Omega_{il}^{\,b}\Omega_{jl}^{\,a}
   -\frac14\sum_c c_{ab}^{\ \ c}\Omega_{ij}^{\,c},\notag\\
 R_{ijab}
 &=\frac{r^2}{4}\sum_l
   (\Omega_{il}^{\,a}\Omega_{jl}^{\,b}
    -\Omega_{jl}^{\,a}\Omega_{il}^{\,b})
   +\frac12\sum_c c_{ab}^{\ \ c}\Omega_{ij}^{\,c},\notag\\
 R_{ijk a}
 &=-\frac r2\big((D_i\Omega)_{jk}^{\,a}-(D_j\Omega)_{ik}^{\,a}\big)
   -\frac r2(\vartheta_i\Omega_{jk}^{\,a}
              -\vartheta_j\Omega_{ik}^{\,a})
   +r\vartheta_k\Omega_{ij}^{\,a},\\
 R_{i a b c}
 &=\frac r2\sum_l\vartheta_l
     (\delta_{ab}\Omega_{il}^{\,c}-\delta_{ac}\Omega_{il}^{\,b}),\notag\\
 R_{abcd}
 &=(r^{-2}-|\vartheta|^2)
       (\delta_{bc}\delta_{ad}-\delta_{ac}\delta_{bd}).\notag
\end{align}
In $R_{iabj}$, the vertical derivative and the vertical Levi--Civita
term contribute respectively $-1/2$ and $+1/4$ times
$\sum_c c_{ab}^{\ \ c}\Omega_{ij}^{\,c}$, giving $-1/4$.
The purely horizontal block is
\begin{align}\label{eq:conn-horizontal-block}
 R_{ijkl}=R^B_{ijkl}+\frac{r^2}{4}\bigl(
 -\langle\Omega_{jk},\Omega_{il}\rangle
 +\langle\Omega_{ik},\Omega_{jl}\rangle
 +2\langle\Omega_{ij},\Omega_{kl}\rangle\bigr).
\end{align}
In particular, for horizontal $X,Y$,
\begin{equation}\label{eq:conn-horizontal-sectional}
 \langle R(X,Y)Y,X\rangle_G
 =\langle R^B(X,Y)Y,X\rangle_{g_B}-\tfrac34r^2|\Omega(X,Y)|^2.
\end{equation}
The remaining components follow from the curvature symmetries.

\subsection{A lower bound on every decomposable bivector}
Let $P=X+U$ and $Q=Y+V$, where $X,Y$ are horizontal for $\omega$ and $U,V$
are components in the unit vertical frame.  In expressions such as
$[U,V]$, identify $\sum_a U_aE_a$ with $\sum_a U_ae_a\in\im\Hh$ and
use the bracket of this fixed Lie algebra.  Decompose
\[
 P\wedge Q=H+B+W,\qquad
 H=X\wedge Y,\quad B=X\otimes V-Y\otimes U,\quad W=U\wedge V,
\]
and put $h=|H|$, $m=|B|$, $k=|W|$.  These are orthogonal components,
so $|P\wedge Q|^2=h^2+m^2+k^2$.  We also regard the mixed tensor $B$
as the linear map from the vertical coefficient space to the horizontal
space given by
\[
 B\xi=X\langle V,\xi\rangle-Y\langle U,\xi\rangle.
\]
Thus $B^*Z=V\langle X,Z\rangle_G-U\langle Y,Z\rangle_G$ for horizontal
$Z$, and $\|B\|_{\mathrm{HS}}=m$, where the Hilbert--Schmidt norm is the
square root of the sum of squared matrix entries in orthonormal bases.
Use the curvature bilinear form
$\mathcal R(A\wedge B,C\wedge D)=\langle R(A,B)D,C\rangle_G$ and put
$\mathcal K(P\wedge Q)=\langle R(P,Q)Q,P\rangle_G$.

The terms linear in $\Omega$ combine as follows.
The mixed block and horizontal--vertical block in
\eqref{eq:connection-curvature-components} give, respectively,
$-\frac12\langle\Omega(X,Y),[U,V]\rangle$ and
$-\langle\Omega(X,Y),[U,V]\rangle$.  Their sum is
\begin{equation}\label{eq:connection-linear-terms}
 -\frac32\langle\Omega(X,Y),[U,V]\rangle.
\end{equation}
For $B_{ia}=X_iV_a-Y_iU_a$, the mixed contribution follows from
\[
 \sum_{i,j,a,b,c}B_{ia}B_{jb}c_{ab}^{\ \ c}\Omega_{ij}^{\,c}
       =2\langle\Omega(X,Y),[U,V]\rangle.
\]
Since $|[U,V]|=2|U\wedge V|$, the absolute value of
\eqref{eq:connection-linear-terms} is at most $3M_0hk$.

To combine the quadratic terms, for a vertical coefficient vector $U$
define the base covector
\[
 \Omega_U(X)(Z)=\langle\Omega(X,Z),U\rangle.
\]
The entire part of the sectional numerator that is quadratic in $\Omega$
is
\begin{equation}\label{eq:connection-quadratic-terms}
 \begin{split}
 \mathcal Q_{\Omega^2}={}&-\frac{3r^2}{4}|\Omega(X,Y)|^2
   +\frac{r^2}{4}|\Omega_V(X)+\Omega_U(Y)|^2\\
 &-r^2\langle \Omega_U(X),\Omega_V(Y)\rangle_{g_B^{-1}}.
 \end{split}
\end{equation}
Indeed, the quadratic mixed block contributes
\[
 \frac{r^2}{4}\bigl(
 |\Omega_V(X)|^2+|\Omega_U(Y)|^2
 -2\langle \Omega_U(X),\Omega_V(Y)\rangle_{g_B^{-1}}\bigr),
\]
and the horizontal--vertical cross block contributes
\[
 \frac{r^2}{2}\bigl(
 \langle \Omega_V(X),\Omega_U(Y)\rangle_{g_B^{-1}}
 -\langle \Omega_U(X),\Omega_V(Y)\rangle_{g_B^{-1}}\bigr).
\]
Their sum, together with \eqref{eq:conn-horizontal-sectional}, proves
\eqref{eq:connection-quadratic-terms}.

For estimates, choose an orthonormal basis $P,Q$ of the two-plane that
diagonalizes the symmetric form
$(Z,Z')\mapsto\langle\pr_{\mathrm H}Z,\pr_{\mathrm H}Z'\rangle_G$
restricted to that plane, where $\pr_{\mathrm H}$ is the horizontal
orthogonal projection.  Then
\begin{equation}\label{eq:conn-principal-pair}
 \langle X,Y\rangle_G=\langle U,V\rangle=0,\qquad
 h=|X||Y|,\quad k=|U||V|,\quad
 m^2=|X|^2|V|^2+|Y|^2|U|^2.
\end{equation}
Under a change of spanning basis, the sectional numerator and each
area squared scale by the square of the determinant.  The resulting
estimates therefore hold for any independent spanning pair.  Since
$|\Omega_U(X)|\leq M_0|X||U|$, discarding the square in
\eqref{eq:connection-quadratic-terms} gives
\begin{equation}\label{eq:conn-quadratic-lower}
 \mathcal Q_{\Omega^2}\geq
 -\tfrac34r^2M_0^2h^2-r^2M_0^2hk.
\end{equation}

To compute the mixed terms with three horizontal arguments, evaluate the
$R_{ijk a}$ block in \eqref{eq:connection-curvature-components} at its repeated horizontal
argument.  For $Z=X$ or $Z=Y$ and any vertical vector $U_0$, alternation
of $\Omega$ and $D\Omega$ in their two-form arguments gives
\[
 \langle R(X,Y)Z,U_0\rangle_G
 =\frac r2\langle(D_Z\Omega)(X,Y),U_0\rangle
  +\frac{3r}{2}\vartheta(Z)\langle\Omega(X,Y),U_0\rangle.
\]
The cross contribution to the curvature quadratic form on $H+B$ is
$2\langle R(X,Y)Y,U\rangle_G-2\langle R(X,Y)X,V\rangle_G$.
Consequently,
\begin{equation}\label{eq:conn-three-horizontal-contraction}
 \begin{split}
 \mathcal Q_{3\mathrm H}={}&r\bigl(
 \langle(D_Y\Omega)(X,Y),U\rangle
 -\langle(D_X\Omega)(X,Y),V\rangle\bigr)\\
 &+3r\bigl(\vartheta(Y)\langle\Omega(X,Y),U\rangle
 -\vartheta(X)\langle\Omega(X,Y),V\rangle\bigr).
 \end{split}
\end{equation}
For the first line, contract the array
$((D_i\Omega)(X,Y)^a)_{ia}$ with $-B_{ia}$; its norm is at most $M_1h$.
For the second, contract $\Omega(X,Y)$ with
$\vartheta(Y)U-\vartheta(X)V=-B^*\vartheta^\sharp$.
Consequently
\[
 |\mathcal Q_{3\mathrm H}|
 \leq r(M_1+3|\vartheta|M_0)hm.
\]
For the mixed terms with one horizontal argument, use $W$ also for the skew
endomorphism of the vertical coefficient space associated with $U\wedge V$:
\[
 W\xi=U\langle V,\xi\rangle-V\langle U,\xi\rangle.
\]
In orthonormal bases put
\[
 W_{ab}=U_aV_b-V_aU_b,\qquad
 Z_{ia}=\sum_l\vartheta_l\Omega_{il}^{\,a}.
\]
The matrices $B,W,Z$ have sizes $d\times3$, $3\times3$, and
$d\times3$, respectively.  Write
$\langle A,B\rangle_{\mathrm{HS}}=\operatorname{tr}(A^*B)$ for the
Hilbert--Schmidt inner product.  The terms with one horizontal argument are
\[
 2\langle R(X,V)V,U\rangle_G-2\langle R(Y,U)V,U\rangle_G.
\]
Substitution of the $R_{iabc}$ block in \eqref{eq:connection-curvature-components}
turns this expression into
\[
 r\sum_{i,a}Z_{ia}\bigl(
 |V|^2X_iU_a+|U|^2Y_iV_a
 -\langle U,V\rangle(X_iV_a+Y_iU_a)\bigr).
\]
The coefficient in parentheses is $-(BW)_{ia}$, so
\begin{equation}\label{eq:conn-one-horizontal-contraction}
 \mathcal Q_{1\mathrm H}=-r\langle BW,Z\rangle_{\mathrm{HS}}.
\end{equation}
If $U,V$ are dependent, $W=0$.  Otherwise $W$ vanishes on their
orthogonal complement and acts on their span as a rotation by a right
angle multiplied by $k$.  Hence $\|W\|_{\mathrm{op}}=k$ in both cases, and
$\|BW\|_{\mathrm{HS}}\leq\|B\|_{\mathrm{HS}}\|W\|_{\mathrm{op}}=mk$.
If $\vartheta=0$, then $Z=0$.  Otherwise choose the first base frame
vector to be $\vartheta^\sharp/|\vartheta|$.  In this frame
$Z_{ia}=|\vartheta|\Omega_{i1}^{\,a}$, so the stated $i<j$ norm gives
\[
 \|Z\|_{\mathrm{HS}}^2
 =|\vartheta|^2\sum_{i=2}^{d}\sum_{a=1}^{3}(\Omega_{i1}^{\,a})^2
 \leq |\vartheta|^2\sum_{i<j,a}(\Omega_{ij}^{\,a})^2
 \leq |\vartheta|^2M_0^2.
\]
Cauchy--Schwarz now gives
$|\mathcal Q_{1\mathrm H}|\leq r|\vartheta|M_0km$.
The Hessian term is $\sum_{i,j,a}\mathsfH_{ij}B_{ia}B_{ja}$.
Here $\lambda_{\min}(\mathsfH)$ denotes the least eigenvalue of the
self-adjoint endomorphism $g_B^{-1}\mathsfH$, so this term is at least
$\lambda_{\min}(\mathsfH)m^2$.  The vertical block is
$(r^{-2}-|\vartheta|^2)k^2$.

Combining these estimates with \eqref{eq:connection-linear-terms} and the
base sectional lower bound gives the estimate
\begin{align}\label{eq:connection-curvature-estimate}
 \mathcal K(P\wedge Q)\geq{}&
 (\kappa_0-\tfrac34r^2M_0^2)h^2+\lambda_{\min}(\mathsfH)m^2
 +(r^{-2}-|\vartheta|^2)k^2\notag\\
 &-r(M_1+3|\vartheta|M_0)hm
 -(3M_0+r^2M_0^2)hk-r|\vartheta|M_0km.
\end{align}
For $C_0=16$, this implies the bound
\begin{equation}\label{eq:connection-lower-estimate}
\begin{aligned}
 \mathcal K(P\wedge Q)\geq{}&
  (\kappa_0-C_0r^2M_0^2)h^2+(r^{-2}-|\vartheta|^2)k^2\\
 &+(\lambda_{\min}(\mathsfH)-C_0r^2M_0^2)m^2\\
 &-C_0r(M_1+|\vartheta|M_0)hm
  -C_0(M_0+r^2M_0^2)hk\\
 &-C_0r|\vartheta|M_0km.
\end{aligned}
\end{equation}
Since $H=X\wedge Y$, the base sectional lower bound gives
$\langle R^B(X,Y)Y,X\rangle_{g_B}\geq\kappa_0h^2$.

\subsection{Estimates for the mixed terms}
Let $C_0=16$, let $\Lambda$ satisfy \eqref{eq:connection-concavity},
and assume $-\Hess\varphi\geq\Lambda g_B$.  Define
\[
 P_* =\max\{1,\|\varphi\|_\infty\},\qquad
 D_* =\max\{1,\|d\varphi\|_\infty\},\qquad
 M_* =\max\{1,M_0\}.
\]
An explicit sufficient scale is
\begin{equation}\label{eq:conn-epsilon-threshold}
 \begin{split}
 \eps_{\mathrm{conn}}=\frac12\min\biggl\{
 &1,\ \frac{\log2}{P_*},\ \frac1{2M_*},\
 \frac1{D_*M_*},\ \frac{\Lambda}{D_*^2},\
 \frac{\kappa_0}{128C_0^2M_*^2},\
 \left(\frac1{2C_0^2D_*^2M_*^2}\right)^{1/3}
 \biggr\}.
 \end{split}
\end{equation}
All denominators are positive, including when $M_0=0$.

For $0<\eps\leq\eps_{\mathrm{conn}}$ set
$r^2=\eps e^{\eps\varphi}$.  Then
\[
 \tfrac12\eps\leq r^2\leq2\eps,\qquad
 \vartheta=\tfrac12\eps\,d\varphi,\qquad
 \mathsfH=-\tfrac12\eps\Hess\varphi
     -\tfrac14\eps^2d\varphi\otimes d\varphi.
\]
In particular, $\lambda_{\min}(\mathsfH)\geq\eps\Lambda/4$.
The threshold also gives
\begin{gather}\label{eq:conn-smallness-consequences}
 r^2M_0\leq1,\qquad |\vartheta|M_0\leq1,\qquad
 C_0r^2M_0^2\leq\kappa_0/4,\notag\\
 r^2|\vartheta|^2\leq1/4,\qquad
 64C_0^2r^2M_0^2\leq\kappa_0,\qquad
 2C_0^2\eps^3\|d\varphi\|_\infty^2M_0^2\leq1.
\end{gather}
Young's inequality gives the following bounds for the three mixed terms
in \eqref{eq:connection-lower-estimate}:
\begin{align}\label{eq:connection-young-estimates}
 C_0r(M_1+|\vartheta|M_0)hm
 &\leq\frac{\kappa_0}{8}h^2
       +4C_0^2\kappa_0^{-1}\eps(M_1+1)^2m^2,\notag\\
 C_0(M_0+r^2M_0^2)hk
 &\leq\frac{\kappa_0}{8}h^2
       +\frac{8C_0^2M_0^2}{\kappa_0}k^2
 \leq\frac{\kappa_0}{8}h^2+\frac{k^2}{8r^2},\\
 C_0r|\vartheta|M_0km
 &\leq\frac{k^2}{8r^2}
       +2C_0^2r^4|\vartheta|^2M_0^2m^2
 \leq\frac{k^2}{8r^2}+\eps m^2.\notag
\end{align}
The first and last estimates use $r^2\leq2\eps$ and
$|\vartheta|\leq\eps\|d\varphi\|_\infty/2$.
The horizontal and vertical coefficients left after these absorptions
are at least $\kappa_0/2$ and $1/(2r^2)$.
The mixed coefficient is at least
\[
 \eps\left(\frac\Lambda4-2C_0M_0^2
       -4C_0^2\kappa_0^{-1}(M_1+1)^2-1\right)
 \geq\frac{\eps\Lambda}{8}
\]
by \eqref{eq:connection-concavity}.  Thus
\begin{equation}\label{eq:conn-final-area-lower}
 \boxed{\displaystyle
 \mathcal K(P\wedge Q)\geq
 \frac{\kappa_0}{2}h^2+
 \frac{\eps\Lambda}{8}m^2+\frac1{2r^2}k^2.}
\end{equation}
Since $r^2\leq2\eps$, division by $h^2+m^2+k^2$ gives
\eqref{eq:connection-margin-main}, including at the boundary.
This proves \cref{thm:connection-metric}.

\section{Choice of the warping functions}\label{app:warping-functions}

We use the polynomial radial reparameterization of
\cite[Section 5.3]{DHZ}.  The estimates below are adapted to our cutoff
interval and to the boundary data in \cref{thm:invariant-transition}.
Fix $a\in(\pi/2,\pi)$, put $F_a=\sin a$, and let
$L_*=\max\{1,L_\rho\}$, where $L_\rho=\sup_y\|K_y\|_{\op}$.
After the southern connection and $A_0>0$ have been fixed, choose a smooth
nondecreasing cutoff $\zeta:[0,\infty)\to[0,1]$ such that
\[
 \zeta=0\text{ on }[0,1],\qquad
 \zeta=1\text{ on }[2,\infty).
\]
Write $M_\zeta=\sup|\zeta'|\geq1$.  Choose $\delta>0$ with
\begin{equation}\label{eq:delta-length}
 \delta<F_a/4
\end{equation}
and
\begin{equation}\label{eq:delta-choice}
 16L_*^2A_0F_a^2\delta^2\leq1,\qquad
 2000L_*^2A_0F_aM_\zeta\delta^3\leq1.
\end{equation}
For example, take
\begin{equation}\label{eq:delta-explicit}
 \delta=\frac12\min\left\{
 \frac{F_a}{4},\ (16L_*^2A_0F_a^2)^{-1/2},
 (2000L_*^2A_0F_aM_\zeta)^{-1/3}\right\}.
\end{equation}

Define $F$ by the strictly increasing polynomial relation
\begin{equation}\label{eq:F-profile}
 s=F+\frac{F^3}{3\delta^2},\qquad F(0)=0,
 \qquad F'=\left(1+\frac{F^2}{\delta^2}\right)^{-1},
\end{equation}
and set
\begin{equation}\label{eq:length-profile}
 \ell=F_a+\frac{F_a^3}{3\delta^2}.
\end{equation}
Then $F(\ell)=F_a$, $F'>0$, $F\leq s$, and
$\ell>F_a>4\delta$.  The inverse of the odd polynomial in
\eqref{eq:F-profile} is odd and smooth, since its derivative is everywhere
positive.  Writing $x=F/\delta$, direct differentiation gives
\begin{equation}\label{eq:F-curvature}
 k_R:=-\frac{F''}{F}=\frac{2}{\delta^2(1+x^2)^3},\qquad
 k_A:=\frac{1-F'^2}{F^2}=\frac{2+x^2}{\delta^2(1+x^2)^2}.
\end{equation}
These formulas include their smooth positive limits at $F=0$.

Use the same southern boundary radius and logarithmic derivative:
\[
 r_a^2=\eps e^{\eps A_0|\cos a|},\qquad
 q_s=\tfrac12\eps A_0F_a,\qquad d=r_aq_s.
\]
Define
\begin{equation}\label{eq:r-profile}
 r(s)=r_a-d\int_s^\ell\zeta(v/\delta)\dd v.
\end{equation}
After $\delta,F,\ell$ have been fixed, let
$\eps_S=\eps_{\mathrm{conn}}$ for the fixed southern data and choose
\begin{equation}\label{eq:epsilon-explicit}
 \eps_* =\frac12\min\left\{\eps_S,
 \frac{\log2}{A_0|\cos a|},
 \frac{1}{A_0F_a\ell},
 \left(\frac{1}{2A_0F_a}\right)^{2/3}\right\},
 \qquad 0<\eps\leq\eps_*.
\end{equation}
The southern metric remains positive.  Also $r_a^2\leq2\eps$,
$q_s\ell\leq1/4$, and $d\leq A_0F_a\eps^{3/2}/\sqrt2<1$.
Consequently
\[
 \tfrac34r_a\leq r\leq r_a,\qquad
 r'=d\zeta(s/\delta)\in[0,1),\qquad
 r''=\frac d\delta\zeta'(s/\delta)\geq0.
\]
In particular $r$ is constant near zero, and
\begin{equation}\label{eq:d-over-r3}
 \frac{d}{r^3}\leq\frac{8q_s}{r_a^2}
 =4A_0F_a e^{-\eps A_0|\cos a|}\leq4A_0F_a.
\end{equation}

For the angular condition in \eqref{eq:north-sufficient}, the ratio of
the term to be subtracted to $k_A$, for $F>0$, satisfies
\begin{equation}\label{eq:polynomial-angular-ratio}
 \frac{2L_\rho^2FF'r'/r^3}{k_A}
 \leq8L_\rho^2A_0F_a\delta^2F\frac{1+x^2}{2+x^2}
 \leq8L_\rho^2A_0F_a^2\delta^2\leq\tfrac12.
\end{equation}
Thus $c_{\mathrm{ang}}\geq k_A/2>0$.  At the center the assertion follows
by continuity, with $r'=0$ nearby.

For the radial condition, $r''$ is supported in $[\delta,2\delta]$.
There $F\leq2\delta$ and hence
\[
 k_R\geq\frac{2}{125\delta^2},\qquad
 \frac{L_\rho^2F^2r''}{r^3}
 \leq16L_\rho^2A_0F_aM_\zeta\delta
 \leq\frac{1}{125\delta^2}\leq\frac{k_R}{2}.
\]
Outside this interval $r''=0$.  Therefore $c_{\mathrm{rad}}\geq k_R/2>0$
on the whole disc.  Since $0\leq F\leq F_a$ and
$(2+x^2)(1+x^2)\geq2$, the two estimates yield the explicit bound
\begin{equation}\label{eq:explicit-north-margin}
 m_N\geq m_N^{\mathrm{poly}}
 :=\frac{1}{\delta^2(1+F_a^2/\delta^2)^3}>0.
\end{equation}
Together with \cref{prop:north-criterion} and $r_{\min}\geq r_a/2$ this
implies
\[
 \Sec_{\bar G_N}\geq
 \frac{1}{\delta^2(1+F_a^2/\delta^2)^3
              (1+4L_\rho^2F_a^2/r_a^2)^2}>0.
\]

Finally,
\[
 F(s)=s-\frac{s^3}{3\delta^2}+O(s^5),\qquad
 \lim_{s\to0}k_R=\lim_{s\to0}k_A=2\delta^{-2}.
\]
Since $r$ is constant near zero, $G_N$ is smooth and has positive
curvature on horizontal planes at the center.  O'Neill's formula then
gives positive curvature of the quotient there.
Since $\ell>4\delta$, the cutoff equals one on a neighborhood of the
northern boundary.  Hence
\[
 F(\ell)=F_a,\qquad r(\ell)=r_a,\qquad r'(\ell)/r_a=q_s.
\]

\section{The Gromoll construction and stabilization}\label{app:parameterized-comparison}

We use the model of framed homeomorphisms in
\cite[Sections 8--9]{SalvatoreTurchin} to compare the Gromoll construction
with stabilization.

\subsection{Families of framed homeomorphisms}
Let $\mathscr H_m=\operatorname{Homeo}_\partial(D^m)$ be the topological
group of homeomorphisms equal to the identity near the boundary,
with the compact-open topology.  We use the $GL_m(\R)$ version
$\mathscr F_m$ of the framed-homeomorphism model in
\cite[Sections 3.4 and 8]{SalvatoreTurchin}.
Its local data are triples $(h,A,H)$ with $h\in\mathscr H_m$, a
continuous map $A:D^m\to GL_m(\R)$, and a continuous family $H_{t,z}$
of invertible germs fixing zero in $\R^m$.  Two representatives of a
germ are identified when they agree on a neighborhood of zero; the
family is understood as a germ along the zero section, so representatives
and their inverses vary continuously in their parameters.  The germ at
$z$ is a map from the normal coordinate over $z$ to that over $h(z)$.
The required endpoints are
\begin{equation}\label{eq:framed-endpoints}
 H(0)=A,\qquad
 H(1)=\Delta h,\qquad
 (\Delta h)_z(v)=h(z+v)-h(z).
\end{equation}
The inclusion $O_m\hookrightarrow GL_m(\R)$ is a homotopy equivalence
and induces the comparison with the $O_m$-framed model.  For triples over $h$ and $k$, composition
has linear part $A_{k(z)}B_z$ and germ path
$H_{t,k(z)}\circ K_{t,z}$, over $h\circ k$.  
All data are the identity on a boundary collar.  The construction may be
made in the simplicial germ model, so the statements below apply to
compact parameter families and their homotopies.

Let $a\in D^r$ be a family parameter.  A relative family
$(h_a,A_a,H_a)$ is extended by the identity in both variables.  Its
underlying product homeomorphism is
\[
 \widehat h(a,z)=(a,h_a(z)).
\]
We associate to $\widehat h$ a triple in $\mathscr F_{m+r}$ with
initial framing $I_r\oplus A_{a,z}$ as follows.  First hold the parameter $a$ fixed and use the path
\[
 (u,v)\longmapsto\bigl(u,(H_a(2\tau))_z(v)\bigr),
 \qquad 0\leq\tau\leq\tfrac12.
\]
Continue with the path
\begin{equation}\label{eq:homeomorphism-path}
 K_t^h(u,v)=\bigl(u,h_{a+tu}(z+v)-h_a(z)\bigr),
 \qquad t=2\tau-1,\quad\tfrac12\leq\tau\leq1.
\end{equation}
At the join $t=0$, the second segment is the direct sum of $I_r$ and
$\Delta h_a$.  At $t=1$, it is $\Delta\widehat h$.
Each $K_t^h$ is an invertible germ: after recovering $u$, its second
coordinate is inverted by
\[
 v=h_{a+tu}^{-1}(h_a(z)+w)-z.
\]
The maps and their inverses depend continuously on all parameters and
are the identity on a common boundary collar.  After rounding the product disc,
this defines $\mathscr P_r:\Omega^r\mathscr F_m\to\mathscr F_{m+r}$.
Forgetting the framing gives the commutative diagram
\begin{equation}\label{eq:framed-forgetful-square}
\begin{tikzcd}[column sep=large]
\Omega^r\mathscr F_m \arrow[r,"\mathscr P_r"] \arrow[d]
 &\mathscr F_{m+r}\arrow[d]\\
\Omega^r\mathscr H_m \arrow[r,"(h_a)\mapsto\widehat h"']
 &\mathscr H_{m+r}.
\end{tikzcd}
\end{equation}
Here $\Omega^r$ is represented by cubical families equal to the identity
near the boundary, with a fixed rounding; this convention is weakly equivalent to the usual based loop
space.  The embedding \eqref{eq:parameter-product-embedding} identifies the
rounded product with the embedded disc containing the support of the
sphere diffeomorphism.  These constructions and their inverses preserve
a common collar for each compact simplex of parameters.

\subsection{Comparison with a smooth Gromoll representative}
For a smooth family $f_a$, the preferred path for a single map is the
difference-quotient path from its derivative to its induced germ
\cite[Lemma 8.1]{SalvatoreTurchin}.  Put $F(a,z)=(a,f_a(z))$.
For $s,t\in[0,1]$, define
\begin{equation}\label{eq:two-parameter-germ}
 \mathcal H_{s,t}(u,v)=
 \begin{cases}
 \displaystyle\left(u,
   \frac{f_{a+stu}(z+tv)-f_a(z)}{t}\right),&t>0,\\[6pt]
 \bigl(u,s(\partial_af_a)_z u+(D_zf_a)v\bigr),&t=0.
 \end{cases}
\end{equation}
For $t>0$, the second component equals
\[
 \int_0^1\!\left[
 s(\partial_a f)_{a+\theta stu,\,z+\theta tv}u
 +(D_zf)_{a+\theta stu,\,z+\theta tv}v
 \right]d\theta.
\]
This expression is smooth in all variables down to $t=0$, where it
has the value prescribed in \eqref{eq:two-parameter-germ}.
For $t>0$, its inverse in the second coordinate is
\[
 v=\frac{f_{a+stu}^{-1}(f_a(z)+tw)-z}{t};
\]
for $t=0$, it is
$(D_zf_a)^{-1}(w-s(\partial_af_a)_z u)$.
For every $t\in[0,1]$, the derivative of $\mathcal H_{s,t}$ at
$(u,v)=(0,0)$ has diagonal blocks $I_r,D_zf_a$ and lower-left block
$s\partial_af_a$, and is therefore invertible.  The parameter-dependent
inverse function theorem supplies a smooth family of inverse germs,
including at $t=0$, on uniform neighborhoods for compact parameter sets.
Its restriction to $t>0$ is the displayed inverse.  For fixed $s$,
follow $t\mapsto\mathcal H_{s,t}$ by the path
$\sigma\mapsto\mathcal H_{s+(1-s)\sigma,1}$, $0\leq\sigma\leq1$.
This gives
\begin{equation}\label{eq:framing-homotopy}
 \Gamma_s(\tau)=
 \begin{cases}
 \mathcal H_{s,2\tau},&0\leq\tau\leq\tfrac12,\\[2pt]
 \mathcal H_{s+(1-s)(2\tau-1),1},&\tfrac12\leq\tau\leq1.
 \end{cases}
\end{equation}
At $\tau=1/2$ both formulas give $\mathcal H_{s,1}$, and
\begin{equation}\label{eq:framing-homotopy-endpoints}
 \Gamma_s(0)=
 \begin{pmatrix}I_r&0\\s\partial_af_a&D_zf_a\end{pmatrix},
 \qquad
 \Gamma_s(1)=\mathcal H_{1,1}=\Delta F
 \quad\text{for every }s.
\end{equation}
At $s=1$, this is the preferred path for $F$, followed by a constant
segment.  At $s=0$, it is the difference-quotient path for $f_a$ with
$a$ fixed, followed by \eqref{eq:homeomorphism-path}.  All paths fix the
product boundary collar.  Reparameterizing both segments by a smooth
function constant near its endpoints makes the concatenation smooth
and preserves its endpoints.
Thus $\Gamma_s$ is a homotopy over $F$ between its associated framed
homeomorphism and $\mathscr P_r$ applied to the framed family of $f_a$.

\begin{samepage}
\begin{lemma}[Compatibility with stabilization]\label{lem:parameterized-comparison}
Let $m=6$ and $r\in\{1,2\}$.  For a smooth $r$-parameter family in $\Diff_\partial(D^m)$ that is the
identity near the entire product boundary, the smoothing comparison of
\cite[Sections 8--9]{SalvatoreTurchin} sends its Gromoll representative
to the stabilization of the family comparison obtained by adding the
identity on $\R^r$.
\end{lemma}
\end{samepage}

\begin{proof}
In \eqref{eq:framed-forgetful-square}, the homotopy
\eqref{eq:framing-homotopy} identifies $\mathscr P_r$ on smooth families with the framed
homeomorphism associated with the Gromoll representative.  Over the
identity family, $h_a=\id$ for all $a$ and $K_t^h=\id$ for every $t$.
Thus $\mathscr P_r$ acts on the fiber by
\[
 (A_{a,z},H_{a,z})\longmapsto
 (I_r\oplus A_{a,z},\ I_r\oplus H_{a,z}),
\]
up to the constant second half of the path.  This takes the direct sum of the linear and microbundle isomorphisms
with the identity on $\R^r$.  After adjunction of the parameter and disc variables it is the
stabilization map of the smoothing homotopy fibers.
The spaces $\mathscr H_{m+r}$ and $\Omega^r\mathscr H_m$ at the bottom
of \eqref{eq:framed-forgetful-square} are contractible: the relative
Alexander contraction contracts $\mathscr H_j$ to the identity, and
its pointwise application contracts the loop space as well.
Use the homotopy-fibration models of
\cite[Sections 8--9]{SalvatoreTurchin}, replacing a projection by its
mapping-path fibration if necessary.  The long exact sequence of a
homotopy fibration over a contractible base makes its fiber inclusion
an isomorphism on all homotopy groups and on path components.
Equivalently, its realization is a weak equivalence; the fiber
inclusion used here is the map $J$ in the proof of Theorem 9.1 of that
paper.  Applying this to the two projections in the square and using
the calculated map on their fibers proves compatibility with
stabilization.

By the fiber sequence used in \cref{lem:Hopf-detection}, the natural
PL/PD-to-TOP comparison induces isomorphisms on $\pi_8$ and $\pi_9$.  Its naturality gives
the same compatibility in the PL model.  After stabilization, the
resulting element of $\pi_j(PL/O)$ is the smoothing class of the sphere
obtained by gluing, under the isomorphism
$\Psi_j:\Theta_j\to\pi_j(PL/O)$ of
\cite[Lemmas 2.5 and 2.7, Proposition 2.10]{CrowleySchick}; see also
\cite[Section 2.4]{Wang}.
\end{proof}

\end{document}